\documentclass[12pt]{amsart}

\usepackage{amsmath, amssymb, amscd}
\usepackage[mathscr]{eucal}
\usepackage{mathrsfs}
\usepackage{verbatim}
\usepackage{version}
\usepackage{nicefrac}
\usepackage{pstricks}
\usepackage{pst-all}
\usepackage{color}
\usepackage[all]{xy}
\usepackage{mathdots}

\usepackage{pict2e}

\usepackage[enableskew]{youngtab}
\usepackage{young}

\usepackage{color}
\usepackage{latexsym} 
\usepackage{amsfonts}
\usepackage{subfigure}
\usepackage{amsthm} 
\usepackage{graphics} 
\usepackage{latexsym}
\usepackage{amscd}
\usepackage{color}
\usepackage{graphicx} 
\usepackage{esint}

\usepackage[colorlinks=true,citecolor=red,linkcolor=blue]{hyper ref}

\usepackage[pagewise, displaymath, mathlines]{lineno}

\newlength\cellsize 
\savebox2{%
\begin{picture}(15,15)
\put(0,0){\line(1,0){15}}
\put(0,0){\line(0,1){15}}
\put(15,0){\line(0,1){15}}
\put(0,15){\line(1,0){15}}
\end{picture}}
\newcommand\cellify[1]{\def\thearg{#1}\def\nothing{}%
\ifx\thearg\nothing
\vrule width0pt height\cellsize depth0pt\else
\hbox to 0pt{\usebox2\hss}\fi%
\vbox to 15\unitlength{
\vss
\hbox to 15\unitlength{\hss$#1$\hss}
\vss}}
\newcommand\tableau[1]{\vtop{\let\\=\cr
\setlength\baselineskip{-16000pt}
\setlength\lineskiplimit{16000pt}
\setlength\lineskip{0pt}
\halign{&\cellify{##}\cr#1\crcr}}}
\savebox3{%
\begin{picture}(15,15)
\put(0,0){\line(1,0){15}}
\put(0,0){\line(0,1){15}}
\put(15,0){\line(0,1){15}}
\put(0,15){\line(1,0){15}}
\end{picture}}
\newcommand\expath[1]{%
\hbox to 0pt{\usebox3\hss}%
\vbox to 15\unitlength{
\vss
\hbox to 15\unitlength{\hss$#1$\hss}
\vss}}

\newenvironment{NB}{
\color{red}{\bf NB}. \footnotesize }{}

\excludeversion{NB} \excludeversion{NB2}

\newcommand{\bF}{\mathbf F}

\newcommand{\bL}{\mathbf L}
\newcommand{\bM}{\mathbf M}
\newcommand{\bN}{\mathbf N}

\newcommand{\cD}{\mathcal D}

\newcommand{\cK}{\mathcal K}

\newcommand{\cM}{\mathcal M}
\newcommand{\cN}{\mathcal N}
\newcommand{\cC}{\mathcal C}
\newcommand{\cO}{\mathcal O}

\newcommand{\cR}{\mathcal R}

\newcommand{\cV}{\mathcal V}
\newcommand{\cX}{\mathcal X}

\newcommand{\fg}{\mathfrak g}

\newcommand{\fo}{\mathfrak o}
\newcommand{\fp}{\mathfrak p}

\newcommand{\ft}{\mathfrak t}

\newcommand{\hR}{\widehat R}

\newcommand{\ocN}{\overline{\cN}}

\newcommand{\rI}{\mathrm{I}}
\newcommand{\rII}{\mathrm{II}}

\newcommand{\sT}{\mathsf T}

\newcommand{\tB}{{\widetilde B}}

\newcommand{\ti}{\widetilde{i}}
\newcommand{\tj}{\widetilde{j}}

\newcommand{\tmu}{\tilde{\mu}}

\newcommand{\tV}{\widetilde V}
\newcommand{\tW}{\widetilde W}

\newcommand{\aaa}{\mathbb A}

\newcommand{\pp}{\mathbb P}

\newcommand{\cc}{\mathbb C}
\newcommand{\rr}{\mathbb R}
\newcommand{\zz}{\mathbb Z}

\newcommand{\ad}{\mathrm{ad} }

\newcommand{\Coeff}{\mathrm{Coeff} }

\newcommand{\dual}{{^\vee} }

\newcommand{\End}{\mathrm{End} }

\newcommand{\gl}{\mathfrak{gl} }

\newcommand{\GL}{\mathrm{GL} }

\newcommand{\git}{/\!\!/ }
\newcommand{\Hom}{\mathrm{Hom} }

\newcommand{\HP}{H\!P}
\newcommand{\Ker}{\mathrm{Ker} }
\newcommand{\Image}{\mathrm{Im} }

\newcommand{\Gr}{\mathrm{Gr} }
\newcommand{\Lie}{\mathrm{Lie} }
\newcommand{\inv}{^{-1}}
\newcommand{\Id}{\mathrm{Id} }
\newcommand{\IH}{I\!H}

\newcommand{\mindeg}{\mathrm{min.deg}\,}
\newcommand{\modality}{\mathrm{mod}}

\newcommand{\rank}{\mathrm{rank} }
\newcommand{\reg}{\mathrm{reg} }

\newcommand{\Res}{\operatornamewithlimits{Res}}

\newcommand{\Spec}{\mathrm{Spec}\, }
\newcommand{\pt}{\mathrm{pt}}

\newcommand{\fsl}{\mathfrak{sl} }

\newcommand{\SL}{\mathrm{SL}}
\newcommand{\SO}{\mathrm{SO}}
\newcommand{\Sp}{\mathrm{Sp}}

\newcommand{\Spin}{\mathrm{Spin}}
\newcommand{\SU}{\mathrm{SU}}

\newcommand{\Sym}{\mathrm{Sym}}
\newcommand{\fsp}{\mathfrak{sp} }

\newcommand{\tr}{\mathrm{tr}}

\newcommand{\USp}{\mathrm{USp} }

\newtheorem{prop}{Proposition}[section]
\newtheorem{thm}[prop]{Theorem}
\newtheorem{lem}[prop]{Lemma}
\newtheorem{cor}[prop]{Corollary}

\theoremstyle{remark}
  \newtheorem{rk}[prop]{Remark}
  
\theoremstyle{definition}
 
 \newtheorem{defn}[prop]{Definition}

\numberwithin{equation}{section}
\numberwithin{figure}{section}

\begin{document}


\title{The $\SO(3)$-instanton moduli space and tensor products of ADHM data}

\author{Jaeyoo Choy}

\address{KIAS, 85 Hoegiro, Dongdaemun-gu, Seoul 02455, South Korea}

\email{choy@kias.re.kr, jaeyoochoy@gmail.com (Primary)}



\begin{abstract}
Let $\cM^K_n$ be the moduli space of framed $K$-instantons with instanton number $n$ over the four-sphere $S^{4}$ when $K$ is a compact simple Lie group of classical type.
Due to Donaldson's theorem \cite{Do}, its scheme structure is given by the regular locus of a GIT quotient of $\mu\inv(0)$ where $\mu$ is the moment map on the associated symplectic vector space of ADHM data.

A main theorem of this paper asserts that $\mu$ is flat for $K=\SO(3,\rr)$ and any $n\ge0$.
Hence we complete the interpretation of the K-theoretic Nekrasov partition function for the classical groups \cite{NS} in term of Hilbert series of the instanton moduli spaces together with the author's previous results \cite{Choy}\cite{Choy2}.

We also write ADHM data for the second symmetric and exterior products of the associated vector bundle of an instanton.
This gives an explicit quiver-theoretic description of the isomorphism $\cM^{K}_{n}\cong\cM^{K'}_{n}$ for all the pairs $K,K'$ with isomorphic Lie algebras.
\end{abstract}

\maketitle
\setcounter{tocdepth}{2}
\tableofcontents

\thispagestyle{empty} \markboth{Jaeyoo Choy}{The $\SO(3)$-instanton moduli space and tensor products of ADHM data}


\section{Introduction}
\label{sec: intro}

The K-theoretic Nekrasov partition function was proposed by Nekrasov \cite{Nek96} to determine the non-perturbative part of the $5$-dimensional $\cN=1,2$ supersymmetric gauge theory. 
For a given gauge group $K$, it is given by the index of coupled Dirac operators on spinors on $K$-instantons (\cite[\S2.1.2]{Nek96}); in particular in algebro-geometric terminology, it is the equivariant integration of the structure sheaf of the moduli scheme in the Grothendieck group of coherent sheaves.
There are several definitions of the moduli space of $K$-instantons as schemes (cf.\ \cite{BE}), but these turn out to be naturally isomorphic to each other as moduli stacks, which makes the K-theoretic Nekrasov partition functions independent from choice of scheme structures of the moduli space. 

Nekrasov and Shadchin \cite{NS} define (a version of) the K-theoretic Nekrasov partition function as an integral formula, i.e., the equivariant integration of some K-theory class defined over the vector space of ADHM quiver representations.
The above K-theory class is set to be the Koszul complex defining the ADHM data in the vector space. 
Hence the Nekrasov-Shadchin integral formula becomes the K-theoretic instanton partition function for the Uhlenbeck partial compactification (Uhlenbeck space, for short) whose scheme structure is given by Donaldson \cite{Do}, provided that the gauge group $K$ is of classical type, but not $\SO(3,\rr),\SO(4,\rr)$.
This fact follows from Crawley-Boevey's result on normality of the moduli space of quiver representations \cite{CB} for $K=\SU(N)$ in a more general context of doubled quiver. 
For $K=\USp(N/2)$ the unitary symplectic group, it is due to Panyushev \cite{Pa} in the context of commuting variety (the case of rank $0$ vector bundles, but the arbitrary rank case easily follows from the base change argument as was mentioned in \cite{ChoyPhD}\cite{Choy}).
For $K=\SO(N,\rr)$, $N\ge5$, this fact is proved by the author \cite{Choy2}.

The
Nekrasov-Shadchin integral formula is also proposed as a deformed Seiberg-Witten prepotential \cite{SW}.
The topological Nekrasov partition function has been understood in this manner, as was established by  Nakajima-Yoshioka \cite{NY} and Nekrasov-Okounkov \cite{NO} independently for $K=\SU(N)$; for an arbitrary compact gauge group by Braverman-Etingof \cite{BE}. 
For further motivations, see \S\ref{subsec: motivation}. 

We have pursued an interpretation of the K-theoretic Nekrasov partition functions as the ``intrinsic'' instanton partition function, more precisely the generating function of Hilbert series for the moduli schemes $\cM^{K}_{n}$.
Our previous results in \cite{Choy}\cite{Choy2} assure that this is indeed the case for $K=\USp(N/2),\ N\in 2\zz_{\ge0}$ or $K=\SO(N,\rr),\ N\ge5$.
The case $K=\SO(4,\rr)$ is rather exceptional because it is not a simple group. 
In this case the Nekrasov-Shadchin integral formula is not the generating function of Hilbert series for $\cM^{K}_{n}$, but that of the symplectic quotient of the space of the ADHM data (ADHM space, for short) (\cite[\S1.3]{Choy2}).
The only remaining case among the classical groups is $K=\SO(3,\rr)$.
Since $\cM^{\SO(3,\rr)}_{n}\cong \cM^{\SU(2)}_{n}$, the instanton partition functions should be identical in principle, but the Nekrasov-Shadchin integral formula are \textit{not} defined intrinsically. 
Indeed in the computation of the integral formula, there is a difference arising from the non-reduced components of the ADHM space for $\SO(3,\rr)$. 
It will be clarified in the next subsection (Theorem \ref{th: main}).

At the point when the Nekrasov-Shadchin formulas for the classical gauge groups admit an algero-geometric interpretation as generating functions of the Hilbert series, we will survey a few interactions with the related topics in \S\ref{subsec: motivation}.


\subsection{Main results}\label{subsec: intro main result}
Let us state the main theorem on $\SO(3,\rr)$-instantons after set-up of the ADHM description of the $\SO(3,\rr)$-instantons (due to Donaldson \cite{Do}).

Let $n$ be an instanton number $n\ge0$.
Let $k:=4n$.
Let $W:=\cc^3,\ V:=\cc^k$ with the standard orthogonal and symplectic forms respectively.
First we consider the space of ordinary ADHM quiver representations
$\bM=\End(V)^{\oplus2}\oplus \Hom(W,V)\oplus \Hom(V,W)$.
The orthogonal and symplectic forms induce the right adjoint maps $\End(V)\to \End(V),\ B\mapsto B^{*}$ and $\Hom(W,V)\to \Hom(V,W),\ i\mapsto i^{*}$.
The $*$-invariant subspace in $\End(V)$ is denoted by $\fp(V)$. 
We define a subspace of $\bM$ as
    $$
    \bN:=\bN_{V,W}:=\{(B_1,B_2,i,j)\in \bM|\, B_1,B_{2}\in \fp(V), j=i^*\}.
    $$
Let $\mu$ be the moment map on $\bN$ with respect to the natural Hamiltonian $\Sp(V)$-action.
Now the $\SO(3,\rr)$-instantons with instanton number $n$ are bijectively corresponds to the stable-costable ADHM quiver representations in $\mu\inv(0)$ via the monad construction.
Here stability (resp.\ costability) means that any $B_{1},B_{2}$-invariant subspace of $V$ containing the image of $i$ coincides with $V$ itself (resp.\ any $B_{1},B_{2}$-invariant subspace of $V$ contained in the kernel of $j$ is trivial). 
The above pair $(\bN,\mu)$ is also defined even if $k=\dim V\in 2\zz$.
We call the ADHM data in $\mu\inv(0)$ as \textit{$\SO(3)$-data}.
This construction works for $\SO(N,\rr)$-instantons (resp.\ $\USp(N/2)$-instantons) by replacing $W=\cc^{3}$ to $\cc^{N}$ (resp.\ further by exchanging the orthogonal/symplectic structures).

A scheme structure of $\cM^{K}_n$ is given by $\mu\inv(0)^\reg/\Sp(V)$ where `$\reg$' denotes the stable-costable locus.
It is the smooth locus of the singular symplectic quotient $\mu\inv(0)\git \Sp(V)$.
Since the K-theoretic partition function is the $\Sp(V)$-invariant part of the equivariant push-forward of the Koszul complex class on $\bN$, it encodes only coherent sheaves over $\bN\git \Sp(V)$, but not over the genuine instanton space $\cM^K_n$.
(Here and hereafter the instanton space abbreviates the instanton moduli space; similarly Gieseker space, monopole space, Nahm space, etc.)
Therefore the K-theoretic partition function coincides with the Hilbert series of $\mu\inv(0)\git \Sp(V)$ if the Koszul complex becomes the free resolution of the structure sheaf $\cO_{\mu\inv(0)}$ as an $\cO_{\bN}$-module.
This amounts to $\mu\inv(0)$ is a complete intersection. 
See the details on the definition of K-theoretic partition function in \cite[\S1]{Choy}.
The following is the first main theorem on the scheme structures of $\mu\inv(0)$ and $\mu\inv(0)\git \Sp(V)$.

\begin{thm}\label{th: main}
Let $\mu$ be the moment map for $\SO(3)$-data with $k=\dim V\in 2\zz_{\ge1}$.
Then we have the followings:

$(1)$
$\mu\inv(0)$ is a non-reduced complete intersection of dimension $\frac{k^2+3k}2$ with the equi-dimensional $\lfloor \frac k4\rfloor+1$ irreducible components.

$(2)$
The GIT quotient $\mu\inv(0)\git \Sp(V)$ is also a non-reduced of dimension $k$ with the equi-dimensional $\lfloor \frac k4\rfloor+1$ irreducible components.
Each irreducible component is birational to  the product of symmetric products $S^{\frac{k-k'}4}\aaa^{4}\times S^{\frac{k'}2}(\aaa^{2}\times \bF)$ where $k'$ runs over the nonnegative integers in the set $k-4\zz_{\ge0}$.
\end{thm}

Here $\lfloor r\rfloor$ denotes the largest integer in $\rr_{\le r}$ where $r\in\rr$ and $\bF$ denotes the fat point $\Spec \cc[x,y,z]/(x,y,z)^{2}$.

The second main result is the quiver description of the scheme-theoretic isomorphisms $\cM^{K}_{n}\stackrel\cong\to \cM^{K'}_{n}$ in terms of ADHM data. 
The precise statement is given in Theorem \ref{th: isom in ADHM}, but it requires a description of a limit of tensor product ADHM data in the ADHM space with respect to a $\cc^{*}$-action and a dilation. 
Here the isomorphisms mean the obvious identifications of instanton spaces for the pairs $(K,K')$ with the same Lie algebras.
There are precisely five pairs $(K,K')$ among the classical groups
	$$
	\begin{aligned}
	&
	(\USp(1),\SU(2)),(\SU(2),\SO(3,\rr)),(\USp(1)\times \USp(1),\SO(4,\rr)),
	\\
	&
	(\USp(2),\SO(5,\rr)),(\SU(4),\SO(6,\rr)).
	\end{aligned}
	$$
For each pair $(K,K')$, the ADHM data of instantons are different.
This implies that the $K$-instanton space has two different scheme structures.
Nevertheless the two different scheme structures are isomorphic as follows:
The natural isomorphism of Lie algebras $\Lie(K)\cong \Lie(K')$ defines an isomorphism between the spaces of associated vector bundles, hence induces a scheme-theoretic isomorphism $\cM^{K}_{n}\cong \cM^{K'}_{n}$ as Gieseker spaces.
The Gieseker space is scheme-theoretically isomorphic to the space of ADHM data via the monad construction.


\subsection{Further motivations}
\label{subsec: motivation}

The $5$-dimensional $\cN=2$ gauge theory has been interacted with several areas of mathematics, in particular, geometry and representation theory.
A perspective from physics \cite{IS}, so-called Intriligator-Seiberg's mirror symmetry in $3$-dimensional $\cN=4$ theory, has motivated mathematical approach to the theory.
It says a certain duality between the two moduli spaces of ALE gravitational instantons and $ADE$ instantons.
The former side is called the Coulomb branch denoted by $\cM_{C}$, while the latter is the Higgs branch denoted by $\cM_{H}$.
In particular the K-theoretic Nekrasov partition function $Z^{K}$ corresponds to $\cM_{C}$ of $3$-dimensional $\cN=4$ theory via the Kaluza-Klein reductions.
This is an argument in physics literatures, e.g.\ Benini-Tachikawa-Xie \cite{BTX}.
Hence in the case $K=\SO(N,\rr)$, $N\in 2\zz$, $Z^{K}$ computes the Higgs branch mirror $\cM_{H}$ of $\cM_{C}$ for the type $D$ quiver gauge theory, as it turns out be the instanton partition function (\cite[Theorem 1.1]{Choy2}).

The above argument is set up as a mathematical conjecture by Hiraku Nakajima \cite{Nak_Coulomb}.  
We briefly outline a part of his formulation narrowly focused on our case, which would hopefully save effort of the reader unfamiliar with this story. 
Our case is studied as one example in \cite[\S3(ii), App.\ A]{Nak_Coulomb}; whilst, his theory itself ranges far wider.
Let $K=\SO(N,\rr)$, $N\in 2\zz$.
Let $(\lambda,\mu)$ be a pair of coweights of $\Spin(N)$, i.e., group homomorphisms $S^{1}\to \Spin(N)$.
This gives a topological type of an $S^{1}$-equivariant $K$-instanton $P$ on $\rr^{4}$ as follows, where $S^{1}$ acts on $\rr^{4}=\cc^{2}$ by $t.(z,w)=(tz,t\inv w)$.
Since there are precisely two $S^{1}$-fixed points $0,\infty$ on the compactification $S^{4}$ of $\rr^{4}$, the fibres $K$ of $P$ at $0,\infty$ are acted by $S^{1}$ via $\lambda,\mu$ respectively.
Now we set $\cM_{C}:=\cM_{C}(\mu,\lambda)$ to be the Uhlenbeck partial compactification of the moduli space of such $P$.
It is known that $\cM_{C}(\mu,\lambda)\neq\emptyset$ if $\lambda\ge\mu$ and $\lambda$ is dominant.

The Higgs branch $\cM_{H}$ is constructed as a quiver variety of Nakajima \cite{NakQ} as follows.
Let $Q=(Q_{0},Q_{1})$ be the quiver of affine $D$ type (the extended Dynkin diagram of $K$), where $Q_{0}, Q_{1}$ denote the sets of vertices and arrows.
Let $\alpha_{i},\Lambda_{i}$, $i\in Q_{0}$ be the simple root and fundamental weight of the untwisted affine Lie algebra associated to $K_{\cc}=\SO(N)$ respectively.
We assign each vertex $i\in Q_{0}$ to the vector spaces $V_{i},W_{i}$ with dimension given by $\lambda=\sum_{i\in Q_{0}}(\dim W_{i})\Lambda_{i}$ and $\mu=\lambda-\sum_{i\in Q_{0}}(\dim V_{i})\alpha_{i}$. 
Here we identified the weights and coweights as $Q$ is simply laced.
Now we set $\cM_{H}:=\cM_{H}(\mu,\lambda)$ the symplectic quotient by $\prod_{i\in Q_{0}}\GL(V_{i})$ of the cotangent space $ T^{*}(\bigoplus_{i\in Q_{0}}\Hom(W_{i},V_{i})\oplus\bigoplus_{a\in Q_{1}} \Hom(V_{t(a)},V_{h(a)}))$, where $t(a),h(a)$ denote the head and tail vertices of an arrow $a$ respectively.
There are various mathematical speculations of the $3$-dimensional mirror symmetry in \cite{Nak_Questions}, among which the following isomorphisms as graded vector spaces are expected: if $\lambda\ge\mu$ and $\mu$ is dominant, $\HP_{0}(\cM_{C})\cong \IH^{*}(\cM_{H})$ and $\HP_{0}(\cM_{H})\cong \IH^{*}(\cM_{C})$, where $\HP_{*},\IH^{*}$ denote the Poisson homology (in Etingof-Schedler \cite{ES}) and the intersection cohomology respectively (\cite[Question 4.1]{Nak_Questions}).

The above proposed Coulomb branch $\cM_{C}$ is expected to be a $K$-monopole space, more precisely a $K$-instanton space over a multi-Taub-NUT space (\cite[\S5(iii)]{Nak_Questions}\cite[\S1(ii)]{BFN}). 
The (1-)Taub-NUT space is $\rr^{4}$ with an ALF metric, so the (Yang-Mills) instantons have a different ADHM-type description, namely the Cherkis bow data \cite[\S2]{Cherkis}. 
It seems that the usual $\SO(N,\rr)$-instanton spaces are one of such monopole spaces as holomorphic symplectic varieties, but in our case $K=\SO(N,\rr)$ we need a further symmetry on the space of the Cherkis bow data, as was informed to us by Nakajima  (cf.\ \cite[A(iii)]{Nak_Coulomb}).

In a general situation of monopole spaces, e.g.\ Cherkis bow varieties (i.e.\ the moduli spaces of bow data), Braverman, Finkelberg and Nakajima propose a definition of $\cM_{C}$ as follows \cite{BFN}. 
Let $\bL$ be a unitary representation of any compact group $K$.
Let $G_{\cO}=G[[t]]$ (the ind-scheme of sections of the trivial $G$-bundle over $\Spec\cc[[t]]$), $G_{\cK}=G((t))$ (the analogous ind-scheme using the Laurent series ring $\cc((t))$ instead of $\cc[[t]]$),  and $\Gr_{G}:=G_{\cK}/G_{\cO}$ (the affine Grassmannian).
Let $\cR$ be the ind-subbundle of $G_{\cK}\times_{G_{\cO}}\bL_{\cO}$ over $\Gr_{G}$ of the elements $[(g,s)]$ satisfying $g.s\in \bL_{\cO}$.
We set $\cM_{C}(G,\bL):=\Spec H^{G_{\cO}}_{*}(\cR)$, where $H^{G_{\cO}}_{*}$ denotes the $G_{\cO}$-equivariant Borel-Moore homology.
Here $H^{G_{\cO}}_{*}(\cR)$ is known to be a commutative algebra with respect to the convolution product via the first projection $p$ and the group-multiplication $m$ in $\Gr_{G}\stackrel{p}\leftarrow G_{\cK}\times_{G_{\cO}} \Gr_{G}\stackrel{m}\to \Gr_{G}$.

A basic motivating example of the above definition of $\cM_{C}$ is a monopole space defined over the ALF space $(-1,1)\times \rr^{3}$, where $(-1,1)$ is the open interval.
It is realized as the moduli space of Nahm data via the Nahm-Hitchin transform (\cite{Nahm}\cite{Hit}\cite{Nak_Nahm}), provided $K$ is of classical type due to Hurtubise-Murray \cite{HM}.
Let us consider the moduli space of $\Lie(K)$-valued Nahm data over $(-1,1)$, where $K$ is any compact group.
According to Bielawski \cite[\S3]{Bie}, it is isomorphic to a sort of commuting varieties in Bezrukavnikov-Finkelberg-Mirkovi\'c \cite{BFM} which comes out from the context of geometric Satake correspondence.
By \cite[Theorem 2.12]{BFM}, this commuting variety, hence the $\Lie(K)$-valued Nahm space, is also isomorphic to $\cM_{C}(\check G,0)=\Spec H^{\check G_{\cO}}_{*}(\Gr_{\check G})$, where $\check G$ denotes the Langlands dual of $G=K_{\cc}$.
In the case of the affine $A_{n-1}$ type quiver $Q=(Q_{0},Q_{1})$, Nakajima and Takayama \cite[\S6.8]{NT} prove that $\cM_{C}(G,\bL)$ is isomorphic to a Cherkis bow variety, where $G=\prod_{i\in Q_{0}}\GL(V_{i})$ and $\bL=\bigoplus_{i\in Q_{0}}\Hom(W_{i},V_{i})\oplus\bigoplus_{a\in Q_{1}} \Hom(V_{t(a)},V_{h(a)})$ as in the above discussion.




\subsection{Flatness of moment map}
\label{subsec: flatness of moment map}

In the following two subsections some technical aspects and strategies in the present paper are briefly introduced.
These are mainly compared to the ideas in the case of higher rank $\SO(N)$, $N\ge4$ in \cite{Choy2}.

In \cite{Choy2} we have proved flatness of the moment map $\mu$ for $\SO(N)$-data, $N\ge4$.
In fact we proved a stronger flatness as we will see below.
First we explain the idea of proof briefly in the case $N=4$ which suffices due to the base change argument.
Note that $\bN$ is a $G$-Hamiltonian vector space of the form $T^{*}X\oplus Y$ where $G=\Sp(V),\ X=\fp(V),\ Y=\Hom(W,V)$ and furthermore that $Y$ is also of cotangent type (i.e.\ $Y=T^{*}\Hom(\cc^{2},V)$).
For $x\in X$, we denote the moment map for $(Y,G^{x})$ by $\mu_{x}$, where $G^{x}$ denotes the stabilizer subgroup of $x$. 
It was observed that $\mu$ is flat if so is $\mu_{x}$ for a generic nilpotent endomorphism in $X$ (cf.\ \cite[(3.5)]{Choy2}).
We proved the latter flatness using the identification of $\fg^{x}:=\Lie(G^{x})$ with the truncated current algebra $\fsl_{2}[z]/(z^{k/2})$ (\cite[Corollary 3.7]{Choy2}).

In the case of the $\SO(3)$-data, the moment map $\mu_{x}$ fails to be flat for a generic endomorphism $x$ and thus the above argument does \textit{not} work any more.
Hence we perform directly an estimate of $\dim \mu\inv(0)$. 
This estimate is rewritten in terms of $\dim \mu_{x}\inv(0)$ due to Lemma \ref{lem: dim}.
By identifying further $Y=\Hom(W,V)=T^{*}V\oplus V$ as symplectic vector spaces, we reduce the estimate of $\dim\mu_{x}\inv(0)$ to that of $\dim\mu_{(x,v)}\inv(0)$ due to the same lemma, where $v\in V$.
In fact the reduced structure of $\mu_{(x,v)}\inv(0)$ is an affine space, which will be explicitly described in Proposition \ref{prop: dim estimate mux0}.

For $N\ge4$, we also proved that $\mu\times E\colon \bN_{V,W}\to \Lie(\Sp(V))\times (\fp(V)\git \Sp(V))$, $X=(B_{1},B_{2},i,j)\mapsto (\mu(X),[B_{1}])$ is flat (\cite[Lemma 2.2]{Choy2}).
(However this is not true for $N=3$ since there is an irreducible component of $\mu\inv(0)$ in Theorem \ref{th: main} which maps to $0$ via $E$.)
By this fact we deduced the factorization property of $\cM^{\SO(N,\rr)}_{n}$ (\cite[Corollary 2.8]{Choy2}).
On the other hand by \cite[Proposition 2.8]{BFM}, if $G$ is simply connected, there is a flat morphism $\cM_{C}(G,\bL)=\Spec H^{G_{\cO}}_{*}(\cR)\to \Spec(H^{*}_{G_{\cO}}(\pt))=\ft/W$ induced by the $H^{*}_{G_{\cO}}(\pt)$-action on $H^{G_{\cO}}_{*}(\Gr_{G})$ (\cite[\S2.6.40]{CG}), where $H^{*}_{G_{\cO}}$, $\ft$ and $W$  denote the $G_{\cO}$-equivariant cohomology, a Cartan subalgebra of $\Lie(G)$ and the Weyl group of $G$ respectively.
The induced map $\overline{E}\colon \cM^{\SO(N,\rr)}_{n}\to \fp(V)\git \Sp(V)$ and the above map $\cM_{C}(G,\bL)\to \ft/W$ have similarity in the sense that they are defined by Poisson-commuting functions (\cite[Theorem 2.15]{BFM}).
Both facts are quite useful in many works, e.g.\ \cite{BFG}\cite{Choy2}\cite{BFN}\cite{NT}, because the flatness gives some \'etale open subsets and geometric properties such as normality can be checked only over these subsets.
We will also use them in the next subsection regarding tensor products.


\subsection{Tensor products of ADHM data}
\label{subsec: tensor}

In \cite{Choy2} the isomorphism $\cM^{K}_{n}\cong \cM^{K'}_{n}$ was constructed in terms of ADHM data for the pair $(K,K')=(\USp(1)\times\USp(1),\SO(4,\rr))$ by giving the ADHM data of tensor products of vector bundles.
Recall from \cite[\S2.8]{Choy2} that the tensor product map at the ADHM data level is only a rational map, which is not defined on the locus of self-tensor products as will be recollected in \S\ref{subsec: tensor prod}.
Via the factorization property the tensor product map is defined in terms of ADHM data over the open subset $(S^{n}\cc\times S^{n'}\cc)^{\circ}$ of $S^{n}\cc\times S^{n'}\cc$ of $n+n'$ points with disjoint support in the commuting diagram  
	$$
	\xymatrix{ 
	\cM^{\SU(N)}_{n}\times \cM^{\SU(N')}_{n'}
	\ar[rrr]^-{\mbox{\tiny tensor product}}  \ar[d]_{\overline{E}\times\overline{E}}
	&&& \cM^{\SU(N+N')}_{nN'+n'N} \ar[d]^{\overline{E}} 
	\\
	S^{n}\cc\times S^{n'}\cc \ar[rrr]^-{\mbox{\tiny weighted sum}}
	&&& S^{nN'+n'N}\cc
	}
	$$
(see \eqref{eq: tensor prod} for the explicit form).
Here the weighted sum maps $(P,P')\in S^{n}\cc\times S^{n'}\cc$ to $N'P+NP'$.
On the other hand, in order to have an ADHM description of the tensor product map for the other pairs $(K,K')\neq(\USp(1)\times\USp(1),\SO(4,\rr))$ with isomorphic Lie algebras, we need to extend it over a larger open subset.

First we illustrate case-by-case how the self-tensor products are involved and then the above enlarged open subset.
The pair $(\USp(1),\SU(2))$ is not the case, as $\cM^{K}_{n}\cong \cM^{K'}_{n}$ is simply induced from the inclusion $\bN_{V,W}\subset \bM_{V,W}$ where $W=\cc^{2}$.
For the other cases $(\SU(2),\SO(3,\rr))$, $(\USp(2),\SO(5,\rr))$ and $(\SU(4),\SO(6,\rr))$, the isomorphisms $\cM^{K}_{n}\cong \cM^{K'}_{n}$ are induced from $F\mapsto S^{2}F$ (the symmetric product) or $F\mapsto \Lambda^{2}F$ (the exterior product).
Thus we need the ADHM datum of the self-tensor product $F^{\otimes 2}$ as a first step because $S^{2}F,\Lambda^{2}F$ are subbundles of $F^{\otimes 2}$.
Now the ADHM data of $S^{2}F,\Lambda^{2}F$ are quiver subrepresentations of the ADHM datum of the self-tensor product $F^{\otimes 2}$ respectively.

The ADHM description of the tensor product map $\cM^{\SU(N)}_{n}\times \cM^{\SU(N)}_{n}\to \cM^{\SU(2N)}_{2nN}$ extends over the open subset of $S^{n}\cc\times S^{n}\cc$ by adding to $(S^{n}\cc\times S^{n}\cc)^{\circ}$ the locus of the diagonals $(P,P)$ such that $P$ has disjoint support.
See \eqref{eq: Phi action} for $n=1$ and Proposition \ref{prop: limit} for $n\ge2$. 
The explicit forms of the ADHM data of $S^{2}F,\Lambda^{2}F$ are given in Theorem \ref{th: isom in ADHM}.



\subsection{Contents}
\label{subsec: intro contents}

This paper is organized as follows.
In \S\ref{sec: Moment maps for SO-data and flatness} we recollect the moment map $\mu$ defining the $\SO(3)$-data and deduce an estimate of fibre dimension of $\mu$.

In \S\ref{sec: flatness of mu for SO3} we prove Theorem \ref{th: main}.
We need to study vector representation of the current algebra $\fsl_{2}[z]$ in \S\S\ref{subsec: vector repn}--\ref{subsec: flatness of mux} mentioned in \S\ref{subsec: flatness of moment map}.
With this preliminary, Theorem \ref{th: main} comes from combination of the above study of vector representation and the fibre dimension estimate in  \S\ref{sec: Moment maps for SO-data and flatness}.

In \S\ref{sec: tensor} we express the ADHM data for tensor products of vector bundles.
First we do for the self-tensor product in \S\ref{subsec: self-tensor} and then for the second symmetric and exterior products in \S\ref{subsec: symm ext ADHM}.
In \S\ref{subsec: construction of isom} the isomorphism $\cM^{K}_{n}\cong \cM^{K'}_{n}$ is written in terms of these ADHM data.

\vskip.3cm \noindent\textit{Acknowledgement.}
The author would like to express the deepest gratitude to Professor Hiraku Nakajima for guiding the author to the instanton spaces and sharing ideas.  
He is grateful to Professors Hoil Kim, Bumsig Kim, Jinsung Park, Hahng-Yun Chu for support and encouragement.


\section{Moment maps for SO-data and the flatness}
\label{sec: Moment maps for SO-data and flatness}

We derive moment maps defining SO-data from the usual moment maps defining ADHM data in \S\ref{subsec: moment maps} and then modify further them in different forms in \S\ref{subsec: another form}.
As a result the moment maps are always canonically defined on symplectic spaces of the form $T^{*}X\times Y$.
In the last subsection \S\ref{subsec: flatness criteria}, we obtain a dimension estimate for the zero fibre of the moment maps.
This gives flatness criteria of the moment maps.

The conventions in the entire part of this paper are as follows.
We are working over $\cc$.
Schemes are of finite type and vector spaces are finite dimensional unless stated otherwise.
Dimension of a scheme of finite type is defined to be the maximum of dimensions of irreducible components.


\subsection{Moment maps for SO-data}
\label{subsec: moment maps}

In this subsection we define moment maps for SO-data.
Some facts will be left as exercises if they come from direct computation, but the details can be found in \cite[\S2]{Choy2}.

Let $(X,\omega)$ be a symplectic manifold with a Hamiltonian $G$-action.
We denote by $\mu_X\colon X\to \fg\dual$ a moment map.
If $X$ is a linear $G$-representation and $\mu_X(0)=0$, $\mu_X$ is given uniquely as $\mu_X(x)=\frac12\omega(\bullet.x,x)$.

Let $V,W$ be vector spaces.
The usual ADHM quiver representations form a symplectic vector space
        $$
        \bM_{V,W}=\End(V)^{\oplus2}\oplus \Hom(W,V)\oplus \Hom(V,W).
        $$
Here the symplectic structures are natural one on the cotangent space $\bM_{V,W}=T^*(\End(V)\oplus \Hom(W,V))$, where $\End(V),\Hom(V,W)$ are identified with their duals via trace respectively.
The moment map with respect to the natural $\GL(V)$-action is given as
    $$
    \mu_{\bM_{V,W}}(B_1,B_2,i,j)=[B_1,B_2]+ij.
    $$
Here we identified the target space $\gl(V)\dual$ of $\mu_{\bM_{V,W}}$ with $\gl(V)$ via trace pairing.

We further assume that there are symplectic and orthogonal forms $(\,,\,)_V,(\,,\,)_W$ on $V,W$ respectively.
For $i\in \Hom(W,V)$ we denote the right adjoint by $i^*\in \Hom(V,W)$ (i.e.\ $(i(w),v)_V=(w,i^*(v))_W$ for $v\in V,w\in W$).
Via $i\mapsto (i,i^*)$, $\Hom(W,V)$ is a symplectic subspace of $T^*\Hom(W,V)=\Hom(W,V)\oplus \Hom(V,W)$ because we can check directly $((i,i^*),(-j^*,j))_{T^*\Hom(W,V)}=2\tr(ij)$ where $j\in \Hom(V,W)$.
This subspace is invariant under the natural $\Sp(V)$-action since $g^*=g\inv$.

For $B\in \End(V)$, $B^*$ denotes the right adjoint.
Let
    $$
    \fp(V):=\left\{B\in\End(V)\middle|B=B^*\right\}
    $$
(the space of symmetric forms).
This is also invariant under the natural $\Sp(V)$-action.
The trace pairing gives an identification $\fp(V)=\fp(V)\dual$.
Note that the Lie algebra $\fsp(V)$ of $\Sp(V)$ is the subspace $\{B\in\End(V)|B=-B^*\}$ (the space of antisymmetric forms).

As a result the vector space  
    $$
    \bN_{V,W}=\fp(V)^{\oplus2}\oplus\Hom(W,V)
    $$
is a $\Sp(V)$-invariant symplectic subspace of $\bM_{V,W}$.
The moment map $\mu_{\bN_{V,W}}$ with respect to $\Sp(V)$ is given as $\mu_{\bN_{V,W}}(B_1,B_2,i)=[B_1,B_2]+ii^*$.
Note that this is the composite of the restriction of $\mu_{\bM_{V,W}}$ and the projection $\gl(V)\dual\to \fsp(V)\dual$.


\subsection{Another form of moment maps}
\label{subsec: another form}

We modify $\mu_{\bN_{V,W}}$ in another forms.
The reason for the modifications is that we need to understand geometry of the zero-fibre of $\mu_{\bN_{V,W}}$.
Recall that $\mu_{\bN_{V,W}}$ is the sum of two moment maps $\mu_{T^{*}\fp(V)}$ and $\mu_{\Hom(W,V)}$ and thus that the zero-fibre is the fibre product of $T^{*}\fp(V)$ and $\Hom(W,V)$ over $\fsp(V)\dual$ via $-\mu_{T^{*}\fp(V)}$ and $\mu_{\Hom(W,V)}$.
Now the standard base change argument can apply once we know the fibres of one of the moment maps explicitly (cf.\ \cite[\S2.4 and Appendix A]{Choy2}).
We will modify the latter moment map.
The computational details in this subsection can be found in \cite[\S2.4]{Choy2}.

The first modification is as follows.
Let us fix any orthogonal basis of $W$ and then identify $\Hom(W,V)=V^{\oplus N}$, where $N=\dim W$.
The moment map $\mu_V$ with respect to $\Sp(V)$ is given by $\mu_V(v)=\frac12(\bullet v,v)$.
And the moment map $\mu_{V^{\oplus N}}$ with respect to $\Sp(V)$ is the $N$-sum of $\mu_V$.

The second modification is as follows.
For simplicity we assume $N=\dim W=3$.
Let $W=W_1\oplus W_2$ be an orthogonal decomposition such that $\dim W_1=2,\dim W_2=1$.
Let $\{e_1,e_2\}$ be a basis of $W_1$ such that $(e_1,e_1)_W=(e_2,e_2)_W=0$.
We take the right adjoints of linear maps in $\Hom(W_{1},V)$ as
    $$
    \begin{aligned}
    *\colon
    \Hom(W_{1},V)=\Hom(\cc e_{1},V)\oplus\Hom(\cc e_{2},V)
    &
    \to \Hom(V,\cc e_{1})\oplus \Hom(V,\cc e_{2}),
    \\
    i=(i_{1},i_{2})
    &
    \mapsto (i_{2}^{*},i_{1}^{*}).
    \end{aligned}
    $$
This induces an isomorphism
	$$
	\Hom(W_{1},V)\cong T^{*}\Hom(\cc e_{1},V)=\Hom(\cc e_{1},V)\oplus\Hom(V,\cc e_{1}),\quad (i_{1},i_{2})\mapsto (i_{1},i_{2}^{*})
	$$
which pull backs the natural symplectic form of the cotangent space $T^{*}\Hom(\cc e_{1},V)$ to the half of the original one of $\Hom(W_{1},V)$.
As a result we identify
	$$\mu_{\Hom(W,V)}=2\mu_{T^{*}V}+\mu_{V}.$$

\begin{rk}\label{rk: modifications of moment maps}

$(1)$
Recall the modification of $\mu_{\Hom(W,V)}$ in the case $N=\dim W\ge4$ in \cite[\S2.4]{Choy2}.
If $N=\dim W=4$, we identified $\mu_{\Hom(W,V)}=2\mu_{T^{*}V^{\oplus2}}$ using a decomposition of $W$ into two complementary maximal isotropic subspaces.
For $N\ge5$ we used an orthogonal decomposition $W=W_{1}\oplus W_{2}$ where $\dim W_{1}=4$ and then the base change argument.

$(2)$
In the case of Sp-data, $V,W$ are orthogonal and symplectic respectively.
By \cite{Pa}, $\mu_{T^{*}\fp(V)}$ is flat and thus we used the base change argument.
\end{rk}


\subsection{Flatness criteria}\label{subsec: flatness criteria}

The moment maps for $\SO$-data and $\Hom(W,V)$ are defined on symplectic spaces of the form $T^{*}X\times Y$.
In this section we give criteria for the moment maps defined on these spaces.
Let $X$ be a $G$-representation and $Y$ be a symplectic smooth scheme with a Hamiltonian $G$-action.
Let $\mu\colon T^{*}X\times Y\to \fg\dual$ be the moment map.
It is the sum of the two moment maps on $T^{*}X$ and $Y$ (denoted by $\mu_{T^{*}X},\mu_{Y}$ respectively).
For any $x\in X$ we denote by $\mu_{x}\colon Y\to (\fg^{x})\dual$ the moment map on $Y$ with respect to $\fg^{x}=\Lie(G^{x})$.
It is given by the composite of the moment map on $Y$ with the natural projection $\fg\dual\to(\fg^{x})\dual$.

\begin{lem}\label{lem: observation}
Let $x\in X$.

$(1)$
The restriction of $\mu_{T^{*}X}$ to $\{x\}\times X\dual$ has the image $(\fg^{x})^{\perp}$ in $\fg\dual$.
In particular its kernel has dimension $\dim X-\dim \fg.x$.

$(2)$
$(\mu_{Y})\inv((\fg^{x})^{\perp})=\mu_{x}\inv(0)$.
\end{lem}

\proof
(1) is immediate from the fact that the dual of the restriction coincides with the $\fg$-action map $\fg\to X,\ \xi\mapsto \xi.x$.

(2) follows from the above definition of $\mu_{x}$ since $(\fg^{x})^{\perp}=\Ker(\fg\dual\to(\fg^{x})\dual)$ .
\qed\vskip.3cm

Let $\cC(G)$ be the set of conjugacy classes $[G']$ of subgroups $G'$ of $G$.
Let $X_{[G']}$ be the locus of $x\in X$ with $G^{x}$ conjugate to $G'$.
It is well-known that there are only finitely many stabilizer subgroups up to conjugacy.
Thus $X$ is decomposed into finitely many strata $X_{[G']}$.

\begin{lem}\label{lem: dim}
Let $S$ be any irreducible $G$-invariant subvariety of $X_{[G']}$ for a given $[G']\in \cC(G)$.
Then we have
    \begin{equation}\label{eq: dim formula}
    \dim\mu\inv(0)\cap(S\times X\dual\times Y)=\dim S-\dim G.x+\dim X+ \dim\mu_{x}\inv(0).
    \end{equation}
where $x$ is any element of $S$.
Hence the dimension of $\mu\inv(0)$ is given as
	$$
	\dim\mu\inv(0)=\max_{[G']\in \cC(G)}\left(\dim X_{[G']}-\dim G.x+\dim X+\dim\mu_{x}\inv(0)\right)
	$$
where $x$ is any element of $X_{[G']}$.
\end{lem}

\proof
Let $d_{x}:=\dim \mu_{x}\inv(0)$.
Recall that $\mu\inv(0)$ is the fibre product of $T^{*}X$ and $Y$ over $\fg\dual$ via $\mu_{T^{*}X}$ and $-\mu_{Y}$.
We restrict this fibre product to $\{x\}\times X\dual\times Y$.
By Lemma \ref{lem: observation} the restriction is identified with the product $\Ker(X\dual\to (\fg_{x})^{\perp})\times \mu_{x}\inv(0)$, proving 
	$$
	\dim\mu\inv(0)\cap(\{x\}\times X\dual\times Y)=\dim X-\dim \fg.x+d_{x}.
	$$ 
By the $G$-equivariancy of $\mu$ we have $\dim\mu\inv(0)\cap(G.x\times X\dual\times Y)=\dim X+d_{x}$.
This proves the dimension formula \eqref{eq: dim formula}.
\qed\vskip.3cm

We give two flatness criteria as a corollary.
Let $X_{(s)}$ be the locus of elements of $X$ with $G$-orbit dimension $s$.

\begin{cor}\label{cor: flatness criterion}
$\mu$ is flat if $\dim X_{(s)}-s+\dim\mu_{x}\inv(0)\le \dim X+\dim Y-\dim G$ for any $x\in X_{(s)}$ and $s\ge0$.

In particular $\mu$ is flat if so is $\mu_{x}$ for any $x\in X$.
\end{cor}

\proof
For the flatness of $\mu$ we need to show $\dim \mu\inv(0)\le 2\dim X+\dim Y-\dim G$, i.e., the central fibre has the dimension less than or equal to the expected one, thanks to the method of associated cones (cf.\ \cite[Appendix D]{ChoyPhD}) and \cite[Chap.\ III, Exer.\ 10.9]{Hart}.
This inequality follows from the above lemma and the obvious inequality
    $$
    \dim X_{[G']}-\dim G.x_{0}+\dim X+\dim\mu_{x_{0}}\inv(0)\le \dim X_{(s)}-s+\dim X+\max_{x\in X_{(s)}}\dim\mu_{x}\inv(0)
    $$
for any $x_{0}\in X_{[G']}$ and any $[G']\in \cC(G)$ such that $X_{[G']}\subset X_{(s)}$.

Let us prove the second statement.
By flatness of $\mu_{x}$ we have $\dim \mu_{x}\inv(0)=\dim Y-\dim \fg^{x}$.
Then the statement follows from plugging this into the inequality of the first criterion.
\qed\vskip.3cm

\begin{defn}\label{def: modality} 
For a $G$-variety $X$ we define \textit{modality} by 
	$$
	\modality(G:X):=\max_{s\ge0}(\dim X_{(s)}-s).
	$$
\end{defn}

\begin{rk}\label{rk: SO4}
The second criterion of the above corollary was used for the flatness of $\mu$ for $\SO(N)$-data, $N\ge4$.
In the case we set $X=\fp(V)$ and $Y=\Hom(W,V)$.
Then $\mu_{x}\colon Y\to (\fg^{x})\dual$ is flat for any generic nilpotent endomorphism $x\in \fp(V)$ (\cite[Corollary 3.7]{Choy2}).
This was argued in terms of modality as follows.
We may assume $N\in 2\zz_{\ge2}$ since the odd rank case follows from the base change.
We write $\Hom(W,V)=T^*\Hom(\cc^r,V)$ where $\cc^r$ is a maximal isotropic subspace of $W$.
Therefore the flatness of $\mu_{x}$ amounts to the inequality $\modality(\fsp(V)^x:\Hom(\cc^r,V))\le \dim \Hom(\cc^r,V)-\dim\fsp(V)^x\ (=(2r-3)k/2)$.
This was shown in \cite[Theorem 3.6]{Choy2}.

However the flatness of $\mu_{x}$ fails for $\SO(3)$-data because $\mu_{x}$ has the central fibre dimension strictly larger than the expected dimension $3k/2$ for any $k\neq0$ as we will see in Theorem \ref{th: dim estimate mu for Hom cc3 V}.

\end{rk}


\section{Flatness of the moment map for $\SO(3)$-data: the proof of Theorem \ref{th: main}}\label{sec: flatness of mu for SO3}

In this section we prove Theorem \ref{th: main}.
Recall that in this theorem $\mu$ is the moment map for $\SO(3)$-data, and $V$ and $W$ are the symplectic and orthogonal vector spaces of dimension $k$ and $N$ respectively.
In the view of the previous subsection, $\mu$ is the moment map on $T^{*}X\times Y$ with respect to the $\Sp(V)$-action, where $X=\fp(V)$ and $Y=\Hom(W,V)$ as in Remark \ref{rk: SO4}.
As we will see, Theorem \ref{th: main} will be deduced from Theorem \ref{th: main 2} which describes geometry of $\mu\inv(0)$ rather than the quotient $\mu\inv(0)\git \Sp(V)$.

To state the latter theorem we need further notation:
$S^k\cc$ denotes the $k^{\mathrm{th}}$ symmetric product of $\cc$.
$\Delta S^k\cc$ denotes the diagonal (i.e.\ the set of unordered $k$-tuples with one-point support).
Let $\eta=(\eta_1,\eta_2,...,\eta_e)$ be a partition of $k$, i.e.\ a decreasing sequence of nonnegative integers with sum $k$.
$S^\eta\cc$ denotes the product $\prod_{n=1}^eS^{\eta_n}\cc$.
$\Delta S^\eta\cc$ denotes $\prod_{n=1}^e\Delta S^{\eta_n}\cc$.
$(S^\eta\cc)_0$ denotes the set of $(z_1,z_2,...,z_e)$ such that the supports of $z_n$ are mutually disjoint.
Let $(\Delta S^\eta\cc)_0:=(\Delta S^\eta\cc)\cap (S^\eta\cc)_0$.

The morphisms $S^\eta\cc\to S^k\cc$ and $\mu\inv(0)\to S^k\cc$ are given respectively by the sum map $(z_1,z_2,...,z_e)\mapsto z_1+z_2+\cdots+z_e$ and $(B_1,B_2,i)\mapsto E(B_1)$ where $E(B_1)$ denotes the unordered $k$-tuple of eigenvalues.

\begin{thm}\label{th: main 2}
$(1)$
$\mu$ is flat.

$(2)$
$\mu\inv(0)$ has the $\lfloor \frac k4\rfloor+1$ irreducible components.
Moreover each irreducible component has an \'etale dense locally closed subset $\mu\inv(0)\times_{S^k\cc}(\Delta S^\eta\cc)_0$ where $\eta$ is a partition of $k$ consisting of only $4$ or $2$.

$(3)$
If $k\in4\zz_{\ge0}$, the regular locus $\mu\inv(0)^{\reg}$ is Zariski dense open in only one irreducible component.
Otherwise $\mu\inv(0)^{\reg}=\emptyset$.
\end{thm}

This section is organized as follows:
In \S\ref{subsec: vector repn} we study the vector representation $\cc^{2}\otimes\cc[z]$ of the current algebra $\fsl_{2}[z]$.
This is preliminary step to estimate $\dim\mu\inv(0)$ in \S\ref{subsec: flatness of mux}.
In \S\ref{subsec: proof of Th 2} we prove Theorem \ref{th: main 2}.
The flatness of $\mu$ will be proven by the dimension estimate of $\mu\inv(0)$.
Combining the flatness with factorization property, we deduce the description of irreducible components.
In \S\ref{subsec: proof of Th} we deduce Theorem \ref{th: main} from Theorem \ref{th: main 2}.
In \S\ref{subsec: Uhlenbeck space} we compare $\mu\inv(0)\git\Sp(V)$ and the Uhlenbeck space via their natural stratifications.  
One can cross-check the description of irreducible components in Theorem \ref{th: main} at set-theoretic level from the stratification.


\subsection{Vector representation of the current algebra $\fsl_{2}[z]$}\label{subsec: vector repn}
As a preliminary step for the flatness of $\mu$ when $\dim W=3$, we need to study vector representation of the current algebra $\fsl_{2}[z]$.
This appeared in \cite[\S3.4]{Choy2} for the case $\dim W=4$.
But there is an essential difference in the case $\dim W=3$: the moment map $\mu_x\colon \Hom(W,V)\to (\fg^x)\dual$ is no more flat (Remark \ref{rk: SO4}).
After the preliminary we will perform a new dimension estimate of the fibres $\mu_x\inv(0)$ in \S\ref{subsec: flatness of mux}.

Let $T:=\cc^{2}=\cc\langle e_{1},e_{2}\rangle$ with the standard symplectic structure.
We denote some tensor products as $\cc$-vector spaces by 
	$$
	T[z]:=T\otimes\cc[z],\ \fsl_{2}[z]:=\fsl(T)\otimes\cc[z],\ \gl_{2}[z]:=\gl(T)\otimes\cc[z].
	$$
Thus their elements have polynomial expression $x=x_{0}+x_{1}z+\cdots+x_{d}z^{d}$ where $x_{i}\in T$ or $\fsl_2$ or $\gl_2$.
Both $T[z],\gl_{2}[z]$ are left $\gl_2[z]$-modules and thus left $\cc[z]$- and $\fsl_2[z]$-modules.
We call $T[z] $ \textit{vector representation} of the current algebras $\fsl_{2}[z]$ and $\gl_{2}[z]$.

We define
    $$
    V_{d}:=T[z]/z^{d+1}T[z],\quad \fg_{d}:=\fsl_2[z]/z^{d+1}\fsl_{2}[z]
    $$
where $d\ge0$.
We can also write them as $V_{d}=T\otimes\cc[z]/(z^{d+1})$ and $\fg_{d}=\fsl_2\otimes\cc[z]/(z^{d+1})$ where $\cc[z]/(z^{d+1})$ is the truncated polynomial algebra.
Note that any $\fsl_2[z]$-orbit in $V_{d}$ coincides with the $\fg_{d}$-orbit since $z^{d+1}$ annihilates any element.

Let $(\,,\,)_{T[z]}$ be the $T[z]$-linear extension of the symplectic form $(\,,\,)_T$.
The induced symplectic structure on $V_d$ is defined by
    $$
    (f,g)_{V_d}:=\Res_{z=0} \frac{(f,g)_{T[z]}}{z^{d+1}}.
    $$
Note that the multiplication map $z$ as an endomorphism of $V_{d}$ is a nilpotent operator with the associated Young diagram $(d+1,d+1)$.
Any generic nilpotent element $B$ of $\fp(V_{d})$ also has the same type Young diagram (see \cite[Lemma C.1]{Choy2} for a general statement).
Therefore by choosing a suitable symplectic basis of $V_{d}$, one can identify $(V,B)=(V_d,z)$, where $k=\dim V=2(d+1)$.
Under this identification we describe $\Sp(V)^{B}$ in the lemma below.
We use the polynomial expression $x=x_{0}+x_{1}z+\cdots+x_{d}z^{d}$ for $x$ in $V_{d}$ or $\gl_{2}\otimes \cc[z]/(z^{d+1})$.
We define \textit{minimal degree}
	$$
	\mindeg x:=\max\left\{n\ge0\middle|z^{d+1-n}x=0\right\}.
	$$
In particular, $x=0$ if and only if $\mindeg x=d+1$.

\begin{lem}\label{lem: Lie stab}
If $B$ is a generic nilpotent element of $\fp(V)$.
Then under the identification $(V,B)=(V_{d},z)$, we have $\Sp(V)^B=G_{d}$, where $G_{d}:=\exp(\fg_{d})$ the adjoint group.
In particular $\Lie(\Sp(V)^{B})=\fg_{d}$.
\end{lem}

\proof
First we notice that $\gl(V_{d})^{z}=\gl_{2}[z]/(z^{d+1})$, where $(z^{d+1})$ is the two-sided ideal of $\gl_{2}[z]$.
In particular any element of $g\in \Sp(V_{d})^{z}$ has a polynomial expression $g=g_{0}+g_{1}z+\cdots+g_{d}z^{d}$.
For any element $g\in \Sp(V_{d})^{z}$ with $g\neq 1$, it suffices to find $\xi\in \fg_{d}$ such that 
	\begin{equation}\label{eq: mindeg exp}
	\mindeg(\exp(-\xi)g-1)>\mindeg (g-1) 
	\end{equation}
due to the induction on the minimal degree of $g-1$.
For, if this claim holds, we can find $\xi,\xi',\xi'',...$ successively, so that $\cdots\exp(-\xi'')\exp(-\xi')\exp(-\xi)g-1$ has minimal degree $d+1$, hence that $g=\exp(\xi)\exp(\xi')\exp(\xi'')\cdots$.

Now we prove the claim.
From $(g(vz^{d}),g(w))_{V_{d}}=(g_{0}(v),g_{0}(w))_{T}$, $v,w\in T$, we have $g_{0}\in \SL_{2}$.
Thus if $g_{0}\neq1$, we may set $\xi\in \fsl_{2}$ with $\exp(\xi)=g_{0}$ as $\exp\colon \fsl_{2}\to \SL_{2}$ is surjective.
Suppose that $g_{0}=1$. 
Let $n:=\mindeg(g-1)$.
Then $1\le n\le d$ and thus $g_{n}\neq0$.
Letting $\xi:=g_{n}z^{n}$, we obtain \eqref{eq: mindeg exp}.
It remains to check $g_{n}\in \fsl_{2}$.
This follows from that for $v,w\in T$, 
	$$
	0=(vz^{d-n},w)_{V_{d}}=(g(vz^{d-n}),g(w))_{V_{d}}=(g_{n}(v),w)_{T}+(v,g_{n}(w))_{T}.
	$$
   
The second assertion on the Lie algebras is immediate from the first one (recall that  it was also shown in \cite[Lemma 3.14]{Choy2} with a similar idea).
\qed\vskip.3cm

We need a description of the stabilizer of a pair $(z,x)$, $x\in V_{d}$.

\begin{prop}\label{prop: gdx}
Let $x\in V_{d}$ and $n:=\mindeg x$.

$(1)$
There exists $\xi\in\fg_{d}^{x}$ with $\mindeg \xi=0$.

$(2)$
For any such $\xi$, $\fg_d^x$ is given as
	$$
	\begin{aligned}
	\fg_{d}^{x} &
	=\cc[z]\xi+z^{d+1-n}\fg_{d}
	\\
	& =
	\cc \xi\oplus\cc \xi z\oplus\cdots\oplus\cc \xi z^{d-n}\oplus\fsl_{2}z^{d+1-n}\oplus \fsl_{2}z^{d+2-n}\oplus\cdots \oplus \fsl_{2}z^{d}.
	\end{aligned}
	 $$

$(3)$
For $s=2(d+1-n)$ $(n=0,1,...,d+1)$, we have
	$$
	(V_{d})_{(s)}
	=
	\left\{x'\in V_{d}\middle|\mindeg x'=n\right\}=
	(T\setminus0)z^{n}\oplus Tz^{n+1}\oplus Tz^{n+2}\oplus\cdots\oplus Tz^{d}.
	$$
For the other values of $s$, $(V_{d})_{(s)}=\emptyset$.

In particular $\modality(\fg_{d}:V_{d})=0$ and it is attained by all the nonempty $(V_{d})_{(s)}$.
\end{prop}

\proof
(1)
We will find $\xi=\xi_{0}+\xi_{1}z+\cdots+\xi_{d}z^{d}\in \fg^x_d$ with $\xi_{0}\neq0$.
We notice first that $\fsl_{2}^{v}$ is $1$-dimensional for any nonzero $v\in\cc^{2}$ because by the $\SL_{2}$-action we may set $v=(1,0)$ and then $\fsl_{2}^{v}$ is spanned by $E$ the elementary matrix with $1$ in the $(1,2)^{\mathrm{th}}$ entry. 
We set $\xi_{0}\in \fsl_{2}$ with $\fsl_{2}^{x_{n}}=\cc\xi_{0}$ (recall $x_{n}\neq0$ as $n$ is the minimal degree of $x$).
We will find $\xi_{1},\xi_{2},...,\xi_{d}\in\fsl_{2}$ inductively.
Assume the induction hypothesis up to $m$, i.e.
    $$
    (\sum_{l=0}^{m}\xi_{l}z^{l})x=0\mod z^{m+n+1}
    $$
as far as $m+n\le d$.
The initial step $m=0$ is already done as $\xi_{0}x\mod z^{n+1}=\xi_{0}x_{n}z^{n}\mod z^{n+1}=0$.
We need to find a solution $\xi_{m+1}$ of the equation
    $$
    \xi_{m+1}x_{n}+\xi_{m}x_{n+1}+\cdots+\xi_{0}x_{n+m+1}=0
    .
    $$
This is always solvable because $x_{n}\neq0$ so that it can be transformed to any vector of $T$ by the $\fsl_{2}$-action.
If $m+n>d$, we set $\xi_{m}$ to be any element of $\fsl_{2}$.

(2)
It is clear that for all $m\ge d+1-n$, $\Coeff_{z^{m}}(\fg_{d}^{x})=\fsl_{2}z^{m}$ because
$\fsl_{2}z^{m}$ itself annihilates $x$.
So it suffices to prove that if $\xi'\in \fg_{d}^{x}$, we have $\xi'\in\cc[z]\xi\mod z^{d+1-n}$.
Let $\xi'\in\fg_{d}^{x}$ with $n':=\mindeg\xi'$.
If $n'\ge d-n+1$, there is nothing to show.
So we suppose $n'\le d-n$.
Since $\xi'_{n'}$ annihilates $x_{n}$, there exists $c\in \cc^{*}$ such that $\xi'_{n'}=c\xi_{0}$.
Here we used $\fsl_{2}^{x_{n}}=\cc\xi_{0}$.
Now $\xi'-c\xi z^{n'}$ has minimal degree $>n'$ and annihilates $x$.
By the induction on minimal degree $n'$, it is contained in $\cc[z]\xi\mod z^{d+1-n}$.
Therefore so is $\xi'$.

(3)
By (2), we have
    $$
    (V_{d})_{(s)}=\left\{x'\in V_{d}\middle|\mindeg x'=n\right\},
    \quad
    s=\dim\fg_{d}-\dim\fg_{d}^{x}=2(d+1-n).
    $$
Note also that $\dim (V_{d})_{(s)}=s$ unless $(V_{d})_{(s)}=\emptyset$.
\qed\vskip.3cm

\begin{rk}
$(1)$
If $n=\mindeg x=0$ in the above proposition, we obtain
    $$
    \fg_{d}^{x}=\cc[z]\xi
    $$
(a principally generated $\cc[z]$-module).
Even for $x\in\Hom(\cc^{r},V_{d})$, $r\ge2$, as far as $\rank x_{0}=1$, $\fg_{d}^{x}$ is also a principally generated $\cc[z]$-module (\cite[Lemma 3.9]{Choy2}).

$(2)$
The modality $0$ in (3) of the proposition can be proven slightly differently: $(V_{d})_{(s)}$ is always one $G_{d}$-orbit, where $G_{d}$ denotes the adjoint group of $\fg_{d}$.
This is immediate from the corollary in the below.
\end{rk}

\begin{cor}\label{cor: orbit}
For the adjoint group $G_{d}=\exp(\fg_{d})$, we have
	$$
	G_{d}.e_{1}z^{n}=\left\{x\in V_{d}\middle|\mindeg x=n\right\}\, (=(V_{d})_{(s)})
	$$
where $n=0,1,...,d+1$ (and $s=2(d+1-n)$).
\end{cor}

\proof
The inclusion $\subset$ is obvious.
To check the opposite we observe that the RHS and any $G_{d}$-orbit in it have the same dimension $2(d+1-n)$ by Proposition \ref{prop: gdx} (2).
Since the RHS is an irreducible Zariski locally closed set, it is one orbit.
\qed\vskip.3cm


\subsection{Dimension estimate of $\mu\inv(0)$ for the moment map $\mu\colon \Hom(\cc^{3},V_{d})\to\fg_{d}\dual$}
\label{subsec: flatness of mux}

By identification $\Hom(\cc^3,V_d)=V_d^{\oplus3}$ by the standard orthogonal basis  of $\cc^3$, $\mu$ is the sum of the three identical moment maps defined on the symplectic vector space $(V_d,(\,,\,)_{V_d})$ with respect to $\fg_d$.
The main theorem in this subsection is the following:

\begin{thm}\label{th: dim estimate mu for Hom cc3 V}
For the moment map $\mu\colon \Hom(\cc^{3},V_{d})\to\fg_{d}\dual$ with respect to the $\fg_d$-action, we have $\dim \mu\inv(0)=4d-2\lfloor \frac d2\rfloor +3$.

In particular if $d=0,1$, $\dim \mu\inv(0)=2\dim V_{d}-1$.
For $d\ge2$, $\dim \mu\inv(0)\le 2\dim V_{d}-3$.
\end{thm}
The proof will appear at the end of this subsection.
We will use Lemma \ref{lem: dim} in the proof, so we need to rewrite $\mu$ as the moment map defined on a space of the form $T^*X\times Y$.
This was done in \S\ref{subsec: another form} where $X=Y=V_{d}$.

To use Lemma \ref{lem: dim} we need dimension estimate of $\mu_x\inv(0)$ for each $x\in V_d$, where $\mu_x\colon V_d\to (\fg_d^x)\dual$ is the  moment map with respect to $\fg_d^x$.
By Corollary \ref{cor: orbit} we may set $x=e_1z^n$.
Let $E,F,H$ be the standard $\fsl_2$-triple:
    $$
    E=\begin{pmatrix} 0 & 1 \\ 0 & 0 \end{pmatrix},\
    F=\begin{pmatrix} 0 & 0 \\ 1 & 0 \end{pmatrix},\
    H=\begin{pmatrix} 1 & 0 \\ 0 & -1 \end{pmatrix}.
    $$

\begin{prop}\label{prop: dim estimate mux0}
Let $x=e_1z^n$, $0\le n\le d+1$.
Then the reduced scheme
    $$
    (\mu_x\inv(0))_{\mathrm{red}}= z^{\lfloor \frac{n-1}2\rfloor+1}\cc[z]e_1+z^{\lfloor \frac d2\rfloor+1}V_{d}.
    $$

In particular if $d=0,1$ and $n=0$, $\mu_{x}\inv(0)$ has codimension $1$ in $V_{d}$.
Otherwise it has codimension $\ge2$.
\end{prop}

\proof
Note first that $\mu_x\inv(0)=\left\{v\in V_d\middle| (v,\fg_d^x.v)_{V_d}=0\right\}$ set-theoretically.
Let $n=0$ first.
By Proposition \ref{prop: gdx} (2), $\fg_d^x=\cc[z]E$ and thus
    $$
    \begin{aligned}
    \mu_x\inv(0)
    &=
    \left\{v\in V_d\middle| (v,E.v)_{V_d}=(v,zE.v)_{V_d}=\cdots=(v,z^dE.v)_{V_d}=0\right\}
    \\
    &=
    \left\{v\in V_d\middle| (v,E.v)_{T[z]}=0\mod z^{d+1}\right\}.
    \end{aligned}
    $$
Using $v=v_0+v_1z+\cdots v_dz^d$, the constraint $(v,E.v)_{T[z]}=0\mod z^{d+1}$ amounts to that for each $0\le l\le d$ the following is zero:
	\begin{equation}\label{eq: Coeff zl}
	\begin{aligned}
	\Coeff_{z^l} & (v,E.v)_{T[z]}
	=
	(v_0,Ev_l)_T+(v_1,Ev_{l-1})_T+\cdots+(v_l,Ev_0)_T
	\\
	&=
	2(v_l,Ev_0)_T+2(v_{l-1},Ev_1)_T+\cdots+\left\{\begin{array}{lll} (v_{l/2},Ev_{l/2})_T && \mbox{$l$: even} \\ 2(v_{(l+1)/2},Ev_{(l-1)/2})_T && \mbox{$l$: odd}\end{array}\right.
	\end{aligned}
	\end{equation}

We claim the solution space for $v$ of this system of quadratic equations is $\cc[z]e_{1}+z^{\lfloor \frac d2\rfloor+1}T[z]\mod z^{d+1}$.
We use the induction on $d$.
When $d=0$, by putting $v_0=ae_1+be_2$ the equation $(v_0,Ev_0)_T=0$ gives $b^2=0$.
So $v_{0}\in \cc e_{1}$.
Let $d\ge1$.
By the induction on $d$ the truncated quadratic equation system 
	\begin{equation}\label{eq: truncated quadratic equation}
	\Coeff_{z^l}(v,E.v)_{T[z]}=0\mod z^{d},\quad 0\le l\le d-1
	\end{equation}
has the solution space $\cc[z]e_{1}+z^{\lfloor \frac {d-1}2\rfloor+1}T[z] \mod z^{d}$.
The original quadratic equation system has one more equation $\Coeff_{z^d}(v,E.v)_{T[z]}=0$ than \eqref{eq: truncated quadratic equation}.
Looking at \eqref{eq: Coeff zl}, when $d$ is odd there is no new equation because $v_{m}\in\cc e_{1}$ for $m\le \frac{d-1}2$ and $Ev_{m}=0$.
Only when $d$ is even,  the equation $(v_{\lfloor\frac d2\rfloor},Ev_{\lfloor\frac d2\rfloor})_{T}=0$  becomes new.
This gives $v_{\lfloor\frac d2\rfloor}\in \cc e_{1}$.
The solution space is now given as in the claim.
This completes the proof of the proposition when $n=0$.

Secondly we assume $n\ge1$.
Then there is the additional defining equation of $\mu_{x}\inv(0)$ 
	\begin{equation}\label{eq: additional defining equation}
	(v,z^{n'}\fg_{d}.v)_{T[z]}=0\mod z^{d+1},\quad n':=d+1-n.
	\end{equation}
This amounts to that \eqref{eq: Coeff zl} vanishes for each $l$ with $E$ replaced by $z^{n'}E,z^{n'}F$ and $z^{n'}H$.
By a similar argument as above, \eqref{eq: additional defining equation} imposes the constraints $v_{l}=0$ for each  $l\ge0$ satisfying $2l+n'\le d$ (i.e.\ $l\le \lfloor\frac{n-1}2\rfloor$).
This finishes the proof of the proposition.
\qed\vskip.3cm

\begin{rk}\label{rk: nonreduced mux}
From the computation $b^2=0$ where $v_0=ae_1+be_2$ in the above proof, we see that $\mu_{x}\inv(0)$ is always a non-reduced scheme.
\end{rk}

Now we can prove Theorem \ref{th: dim estimate mu for Hom cc3 V}.
\vskip.3cm

\textit{Proof of Theorem \ref{th: dim estimate mu for Hom cc3 V}.}
By Lemma \ref{lem: dim} we need to compute $\dim X_{[G']}-\dim G.x$ ($x\in X_{[G']}$) and $\dim \mu_{x}\inv(0)$ where $[G']\in\cC(G)$ and $x\in X_{[G']}$.
By Proposition \ref{prop: gdx} (3), $\dim X_{[G']}-\dim G.x=0$ for any $[G']\in \cC(G)$ and $x\in X_{[G']}$.
By Proposition \ref{prop: dim estimate mux0},  $\dim \mu_{x}\inv(0)$ attains the maximum $2(d-\lfloor \frac d2\rfloor)+1$ precisely when $\mindeg x=0$.
\begin{NB} Simplify $\dim \mu_{x}\inv(0)=\lfloor \frac d2\rfloor+1+2(d-\lfloor \frac d2\rfloor)$ where the RHS is obtained from $n=0$ in the proposition. \end{NB}
Therefore by the dimension formula in Lemma \ref{lem: dim} we obtain $\dim \mu\inv(0)=4d-2\lfloor \frac d2\rfloor +3$.
\begin{NB} It comes from $\dim \mu\inv(0)=\dim V_{d}+2(d-\lfloor \frac d2\rfloor)+3=2(d+1)+2(d-\lfloor \frac d2\rfloor)+1.$\end{NB}
\qed\vskip.3cm


\subsection{Proof of Theorem \ref{th: main 2}}\label{subsec: proof of Th 2}
We get back to the setting of Theorem \ref{th: main 2}: $\mu\colon T^{*}X\times Y\to \fsp(V)\dual$ is the moment map with respect to $\Sp(V)$ where $X=\fp(V),\ Y=\Hom(W,V)$ and $W=\cc^{3}$ (the orthogonal vector space).

First we prove $\mu$ is flat.
We need the factorization property.
To emphasize $k=\dim V$ we use the notation $\mu_{(k)}$ instead of $\mu$.

\begin{lem}\label{lem: factorization} $($cf.\ \cite[Lemma 2.7 (2)]{Choy2}$)$
Let $\eta=(\eta_1,\eta_2,...,\eta_e)$ be a partition of $k$.
Then there is a surjective smooth morphism
    $$
    \sigma\colon \Sp(V)\times\left(\mu_{(\eta_1)}\inv(0)\times \mu_{(\eta_2)}\inv(0)\times\cdots\times \mu_{(\eta_e)}\inv(0)\right)\times_{S^\eta\cc}(S^\eta\cc)_0\to \mu_{(k)}\inv(0)\times_{S^k\cc}(S^\eta\cc)_0
    $$
with the (scheme-theoretic) fibres
    $$
    \sigma\inv(\sigma(g,x))=\left\{(gh\inv,h.x)\middle|h\in \Sp(\eta_1/2)\times\Sp(\eta_2/2)\times\cdots\times\Sp(\eta_e/2)\right\}.
    $$

Moreover $\sigma$ satisfies $\sigma(g,x)=g.\sigma(e,x)$ where $e$ denotes the identity of $\Sp(V)$ and the $\Sp(V)$-action is trivial on the factor $S^{\eta}\cc$.
\end{lem}

This lemma holds similarly for the ordinary ADHM data by replacing $\Sp(V)$ into $\GL(V)$.
But in our case, only when $\eta_1,\eta_2,...,\eta_e$ are even, $\sigma$ has the nonempty domain and target spaces.

If $e\ge2$ and $\eta_e\neq0$, this factorization property gives the induction hypothesis on $k=\dim V$ assuring the \'etale open subset associated to $\eta$, $\mu_{(k)}\inv(0)\times_{S^k\cc}(S^\eta\cc)_0$ has the expected dimension $\dim\fp(V)+2\dim V$.
Thus to prove flatness of $\mu_{(k)}$ we need only the dimension estimate when $e=1$, i.e.,
    $$
    \dim\mu\inv(0)\cap(S\times X\times Y)\le \dim\fp(V)+2\dim V,
    $$
where
    $$
    S:=\left\{B\in\fp(V)\middle| \mbox{$B$ has only one eigenvalue}\right\}.
    $$
By the dimension formula \eqref{eq: dim formula}
this amounts to
    \begin{equation}\label{eq: modified dim formula}
    \dim S_{[G']}-\dim \Sp(V).x+\dim\mu_{x}\inv(0)\le 2\dim V
    \end{equation}
 where $[G']\in \cC(G)$ and $x\in S_{[G']}$.

Let us estimate the LHS of \eqref{eq: modified dim formula}.
Let $x\in S$.
Thus $x$ has only one eigenvalue, say $a$.
It is well-known that there is a $x$-stable decomposition into symplectic subspaces $V=V_{1}\oplus V_{2}\oplus\cdots\oplus V_{l_{0}}$ such that the nilpotent endomorphism $(x-a)|_{V_{l}}$ corresponds to the partition $(\frac{\dim V_{l}}2,\frac{\dim V_{l}}2)$ for each $l$ (\cite[Corollary C.3]{Choy2}).
Note that the conjugacy class of stabilizers $G'=\Sp(V)^x$ is uniquely given by the unordered $l_{0}$-tuple of such partitions.
This implies $\Sp(V)$-orbits in $S_{[G']}$ are in one-to-one correspondence with the eigenvalues $\cc$, which means 
    \begin{equation}\label{eq: dim S minus orbit}
    \dim S_{[G']}-\dim \Sp(V).x=1
    \end{equation}
for any $x\in S_{[G']}$.

On the other hand $\dim \mu_x\inv(0)$ is estimated as follows.
Let $x_{l}:=x|_{V_{l}}\in \fp(V_{l})$.
Since the product $\prod_{l=1}^{l_{0}}\Sp(V_{l})^{x_{l}}$ is a subgroup of $G'$, $\mu_{x}\inv(0)$ is contained in the product $\prod_{l=1}^{l_{0}}\mu_{x_{l}}\inv(0)$ where $\mu_{x_{l}}\colon \Hom(W,V_{l})\to(\fsp(V_{l})^{x_{l}})\dual$ is the moment map with respect to $\Sp(V_{l})^{x_{l}}$.
Since we have a natural identification $\fsp(V_{l})^{x_{l}}=\fg_{d'}$ where $d'=\dim V_{l}/2-1$ by Lemma \ref{lem: Lie stab}, we have the dimension estimate of $\mu_{x_{l}}\inv(0)$ for each $l$ by Theorem \ref{th: dim estimate mu for Hom cc3 V}. 
To sum up we obtain
    \begin{equation}\label{eq: further mod dim formula}
    \dim \mu_{x}\inv(0)\le \sum_{l=1}^{l_{0}}\dim \mu_{x_{l}}\inv(0)\le 2\dim V-l_{0}.
    \end{equation}
By adding \eqref{eq: dim S minus orbit} and \eqref{eq: further mod dim formula}, we obtain an estimate of the LHS of \eqref{eq: modified dim formula}: $\mathrm{LHS}\le 2\dim V-l_0+1$.
So \eqref{eq: modified dim formula} is checked.
The proof of Theorem \ref{th: main 2} (1) is done.

We prove Theorem \ref{th: main 2} (2).
We need to check
	$$
	\begin{aligned}
	&
	\dim\mu\inv(0)\times_{S^k\cc}(\Delta S^\eta\cc)_0=\dim\mu\inv(0)\ \mbox{if}\ \eta_1\le 4,
	\\
	&
	\dim\mu\inv(0)\times_{S^k\cc}(\Delta S^\eta\cc)_0<\dim\mu\inv(0)\ \mbox{if}\ \eta_1\ge6.
	\end{aligned}
	$$ 
As in the proof of (1) we use the dimension formula \eqref{eq: dim formula}.
By the second statement of Theorem \ref{th: dim estimate mu for Hom cc3 V}, \eqref{eq: further mod dim formula} is further modified to
    \begin{equation}\label{eq: further further mod dim formula}
    \dim \mu_{x}\inv(0)\le 2\dim V-\#\left\{1\le l\le l_{0}\middle| \dim V_l=\mbox{$2$ or $4$}\right\}-3\#\left\{1\le l\le l_{0}\middle| \dim V_l>4\right\}.
    \end{equation}
Therefore we have 
	$$
	\begin{aligned}
	&
	\dim \mu_{x}\inv(0)= 2\dim V-1\ \mbox{if $l_{0}=1$ and $\dim V\le 4$},
	\\
	&
	\dim \mu_{x}\inv(0)\le 2\dim V-2\ \mbox{otherwise.}
	\end{aligned}
	$$ 
By \eqref{eq: dim formula} combined with the factorization property (Lemma \ref{lem: factorization}), the \'etale open subset $\mu\inv(0)\times_{S^k\cc}(S^\eta\cc)_0$ has dimension equal to $\dim\mu\inv(0)$ (resp.\ strictly less than $\dim \mu\inv(0)$) if $\eta_1\le4$ (resp.\ otherwise).

It remains to show $\mu_{(4)}\inv(0)$ has two irreducible \'etale locally closed subsets $\mu_{(4)}\inv(0)\times_{S^4\cc}(S^{(2,2)}\cc)_0$ and $\mu_{(4)}\inv(0)\times_{S^4\cc}\Delta S^{(4)}\cc$.
This is equivalent to the following:

\begin{lem}\label{lem: N3k4} $($\cite[Corollary 8.9]{Choy}$)$
There are precisely two irreducible components of $\mu_{(4)}\inv(0)$.
They are set-theoretically the $(\aaa^{2}\rtimes\SL_{2})\times\Sp(2)\times\SO(3)$-orbit closures of some elements $(B_{1}^{(\rI)},B_{2}^{(\rI)},i^{(\rI)}),(B_{1}^{(\rII)},B_{2}^{(\rII)},i^{(\rII)})$ such that $B_{1}^{(\rI)}$ has only one eigenvalue while $B_{1}^{(\rII)}$ has two distinct eigenvalues.
\end{lem}
Here $\aaa^{2}\rtimes \SL_{2}$ acts on $\fp(V)^{\oplus2}$ by
	$$
	\left((a_{1},a_{2}),\begin{pmatrix} a & b \\ c & d\end{pmatrix}\right).(B_{1},B_{2})=(a_{1},a_{2})+(aB_{1}+bB_{2},cB_{1}+dB_{2}).
	$$
The explicit forms of $(B_{1}^{(\rI)},B_{2}^{(\rI)},i^{(\rI)}),(B_{1}^{(\rII)},B_{2}^{(\rII)},i^{(\rII)})$ are given in \cite[Lemma 8.8 and pp.298--299]{Choy}.
This completes the proof of (2).

We prove Theorem \ref{th: main 2} (3).
By \cite[Lemma 8.10 (2)]{Choy},  $(B_{1}^{(\rII)},B_{2}^{(\rII)},i^{(\rII)})$ is a regular element.
By \cite[Lemma 2.5]{Choy2}, any tuple of regular elements maps to a regular element in $\mu\inv(0)$ via $\sigma$ in Lemma \ref{lem: factorization}.
Therefore if $k\in 4\zz_{\ge0}$, the product of $\mu_{(4)}\inv(0)^{\reg}$ maps via $\sigma$ into a Zariski dense open subset of an irreducible component of $\mu\inv(0)^{\reg}$.
To complete the proof of (3) we observe that any element of $\mu_{(2)}\inv(0)$ has nontrivial stabilizer in $\Sp(1)$.
For, if $(B_{1},B_{2},i)\in \mu_{(2)}\inv(0)$, both $B_{1},B_{2}$ are scalars and $i$ is of rank $1$ (see \cite[Corrigendum and addendum:\ Remark 2.4]{Choy}).
Hence for any partition $\eta$ with $\eta_{e}=2$ and $x\in\mu_{(\eta_1)}\inv(0)\times \mu_{(\eta_2)}\inv(0)\times\cdots\times \mu_{(\eta_e)}\inv(0)$, $\sigma(e,x)$ has nontrivial stabilizer in $\Sp(V)$ by the second statement of Lemma \ref{lem: factorization}.
Hence it cannot be a regular element.
This finishes the proof of Theorem \ref{th: main 2} (3).


\subsection{Proof of Theorem \ref{th: main}}\label{subsec: proof of Th}
The first item (1) of the theorem, except the non-reducedness is immediate from the flatness of $\mu$ as $\mu$ has the equi-dimensional fibres.
The non-reducedness comes from the binational description of the irreducible components in (2).

We prove (2).
Note first that the sum map $(\Delta S^{\eta}\cc)_{0}\to S^{k}\cc$ is the $\zz_{n_{1}}\times\zz_{n_{2}}$-quotient map onto the image where $\eta=(4^{n_1},2^{n_2}),\  4n_{1}+2n_{2}=k$.
By Theorem \ref{th: main 2} (2), each irreducible component of $\mu\inv(0)$ contains the $\zz_{n_{1}}\times\zz_{n_{2}}$-quotient of $\mu\inv(0)\times_{S^{k}\cc}(\Delta S^{\eta}\cc)_{0}$ as a Zariski dense locally closed subset.

To describe the above Zariski dense locally closed subset, we consider first $\mu\inv(0)\times_{S^{k}\cc}(S^{\eta}\cc)_{0}$.
By the factorization property (Lemma \ref{lem: factorization}), it is isomorphic to an affine open subset of the free GIT quotient
	\begin{equation}\label{eq: free GIT quotient}
	\left(\Sp(V)\times \mu_{(4)}\inv(0)^{n_1} \times\mu_{(2)}\inv(0)^{n_2}\right)/\Sp(2)^{n_{1}}\times\Sp(1)^{n_{2}}.
	\end{equation}
Recall the quotient is defined with respect to the $\Sp(2)^{n_{1}}\times\Sp(1)^{n_{2}}$-action $h.(g,x)=(gh\inv,h.x)$.
Since the actions of $\Sp(V)$ and $\Sp(2)^{n_{1}}\times\Sp(1)^{n_{2}}$ on $\Sp(V)\times \mu_{(4)}\inv(0)^{n_1} \times\mu_{(2)}\inv(0)^{n_2}$ commute, the GIT $\Sp(V)$-quotient of the above affine scheme \eqref{eq: free GIT quotient} is written as   
	\begin{equation}\label{eq: GIT quotient}
	\left(\mu_{(4)}\inv(0)\git\Sp(2)\right)^{n_1}
	\times
	\left(\mu_{(2)}\inv(0)\git\Sp(1)\right)^{n_2}.
	\end{equation}
We recollect from \cite{Choy} the scheme structures of the factors $\mu_{(4)}\inv(0)\git\Sp(2)$ and $\mu_{(2)}\inv(0)\git\Sp(1)$.
Let $\rho\colon \Hom(W,V)\to \fsp(V),\ i\mapsto ii^*$ where $\dim V=2$.
Let $\cN$ be the minimal nilpotent $\SO(3)$-orbit in $\fo(3)$ and $\ocN$ be its Zariski closure.

\begin{lem}\label{lem: 2 4}
There are (natural) isomorphisms 
	$$ 
	\begin{aligned}
	&
	\mu_{(2)}\inv(0)\git\Sp(1)\cong\aaa^2\times \rho\inv(0)\git\Sp(1)\cong \aaa^{2}\times \bF,
	\\
	&
	\quad \cM^{\SO(3,\rr)}_{1}\cong \aaa^{2}\times \cN\cong \aaa^{2}\times (\aaa^{2}\setminus0).
	\end{aligned}
	$$ 
\end{lem}

\proof
The first isomorphism is given in \cite[Corrigendum and addendum:\ \S3.1]{Choy}.
The second one will be proven in Corollary \ref{cor: Hilbert series}.
The third one is given in \cite[Lemma 8.13]{Choy}.
The last one comes from $\ocN\cong \Spec \cc[x^{2},xy,y^{2}]$ (\cite[Lemma 8.11]{Choy}).
\qed\vskip.3cm

As a result  the affine GIT $\Sp(V)$-quotient $\mu\inv(0)\times_{S^{k}\cc}(\Delta S^{\eta}\cc)_{0}\git \Sp(V)$ is birational to $(\aaa^{4})^{n_{1}}\times (\aaa^{2}\times \bF)^{n_{2}}$.
Now the $\zz_{n_{1}}\times \zz_{n_{2}}$-quotient $S^{n_{1}}\aaa^{4}\times S^{n_{2}}(\aaa^{2}\times \bF)$ is birational to the irreducible component of $\mu\inv(0)\git \Sp(V)$ indexed by $\eta=(4^{n_{1}},2^{n_{2}})$.
This completes the proof of Theorem \ref{th: main}.


\subsection{Uhlenbeck space}\label{subsec: Uhlenbeck space}

Let us consider the Uhlenbeck space of $\cM^{\SO(3,\rr)}_{k/4}$ in the case $k\in 4\zz$.
It has the stratification 
	$$
 	\bigsqcup_{0\le k'\le k,\, k'\in 4\zz} \mu_{(k')}\inv(0)^\reg/\Sp(V_{k'})\times S^{\frac{k-k'}4}\aaa^2 .
 	$$

On the other hand by identification of closed $\Sp(V)$-orbits we have stratification (cf.\ \cite[Theorem 2.6 (1)]{Choy}):
    	\begin{equation}\label{eq: stratification}
    	\mu\inv(0)\git G
	=\bigsqcup_{0\le k'\le k,\, k'\in 4\zz} \mu_{(k')}\inv(0)^\reg/\Sp(V_{k'})\times S^{\frac{k-k'}2}\aaa^2.
    	\end{equation}
Thus $\mu\inv(0)\git G$ is different from the Uhlenbeck space unlike the case $\SO(N)$, $N\ge5$.

From \eqref{eq: stratification} one can check that each stratum indexed by $k'$ is birational to the product of symmetric products $S^{\frac k4}\aaa^{4} \times S^{\frac{k-k'}2}\aaa^{2}$.
This reconfirms Theorem \ref{th: main} (2) when we consider only the reduced scheme structures. 
The case when $k\in 2\zz\setminus 4\zz$ is similar.


\section{Tensor products of ADHM data}
\label{sec: tensor}

The following is the complete list of pairs of simple compact classical groups having the isomorphic Lie algebras:
    $$
    (K,K')= (\USp(1),\SU(2)),(\SU(2),\SO(3,\rr)),(\USp(2),\SO(5,\rr)),(\SU(4),\SO(6,\rr)).
    $$
Therefore there are isomorphisms $\cM^K_n\cong\cM^{K'}_n$ mapping the associated vector bundles
    $$
    F\mapsto F,\ad F,(\Lambda^2 F)_0,\Lambda^2 F
    $$
respectively.
Here $\ad F$ is the trace-free part of $\End(F)$, and $(\Lambda^2 F)_0$
is the kernel of the natural symplectic form $\Lambda^2F\to \cO$.
For a rank $4$ vector bundle $F$ with $\cO\cong \Lambda^{4}F$, the natural orthogonal structure on $\Lambda^{2}{F}$ is given by wedge product.

In this section we interpret the above isomorphisms as the morphisms between the spaces of ADHM data (Theorem \ref{th: isom in ADHM}).
Except the obvious case $(\USp(1),\SU(2))$ the isomorphisms involve self-tensor products of vector bundles $F^{\otimes2}$.
The ADHM datum of tensor product $F\otimes F'$ in terms of ADHM data of $F,F'$ is called \textit{tensor product of ADHM data}.
But we hope the readers do not confuse this with the usual tensor product as representations of the quiver algebra with respect to the diagonal action.

In \S\ref{subsec: tensor prod} we recollect the construction of tensor product in \cite[\S2.8]{Choy2}.
In \S\ref{subsec: self-tensor} we construct the ADHM data of the self-tensor product.
In \S\ref{subsec: symm ext ADHM}, using the self-tensor product ADHM data, we construct the symmetric product and exterior product of ADHM data.
In \S\ref{subsec: construction of isom} we construct the above claimed morphisms in terms of ADHM data.


\subsection{Recollection on tensor product of ADHM data}\label{subsec: tensor prod}

First we recollect tensor product of ADHM data in \cite[\S2.8]{Choy2}.
This was considered for the semisimple classical group pair $(\USp(1)\times\USp(1),\SO(4,\rr))$.

The prototype construction is for the ordinary ADHM data.
We are given two pairs of vector spaces $(V,W)$ and $(V',W')$ as the representation spaces of the ADHM quiver algebra.
Let
	$$
	\tV:=V\otimes W' \oplus W\otimes V',\quad \tW:=W\otimes W' .
	$$
We define a rational map
	\begin{equation}\label{eq: tensor prod}
	\begin{aligned}
	&
	\sT\colon\bM_{V,W}\times\bM_{V',W'}\dasharrow \bM_{\tV,\tW},
	\\
	&
	(x,x')=((B_{1},B_{2},i,j),(B_{1}',B_{2}',i',j'))\mapsto
	\\
	&
	{\small
	\left(
	\begin{pmatrix} B_{1}\otimes \Id_{W'} & 0 \\ 0 & \Id_{W}\otimes B_{1}'\end{pmatrix},
	\begin{pmatrix} B_{2}\otimes \Id_{W'} & \tB^{(1,2)} \\ \tB^{(2,1)} & \Id_{W}\otimes B_{2}' 	 	 \end{pmatrix},
	\begin{pmatrix} i\otimes\Id_{W'} \\ \Id_{W}\otimes i' \end{pmatrix},
	\begin{pmatrix} j\otimes \Id_{W'} & \Id_{W}\otimes j'\end{pmatrix}
	\right)
	}
	\end{aligned}
	\end{equation}
where $\tB^{(1,2)},\tB^{(2,1)}$ satisfy the ADHM equations
	\begin{equation}\label{eq: matrix eq}
	\begin{aligned}
	&
	(B_{1}\otimes \Id_{W}) \tB^{(1,2)}-\tB^{(1,2)} (\Id_{W}\otimes B_{1}')+i\otimes j'=0,
	\\
	&
	(\Id_{W}\otimes B_{1}') \tB^{(2,1)} -\tB^{(2,1)} (B_{1}\otimes \Id_{W'})+j\otimes i'=0.
	\end{aligned}
	\end{equation}
Looking at the first equation, $\tB^{(1,2)}$ is uniquely defined if and only if $B_{1},B_{1}'$ do not have a common eigenvalue.
Its matrix elements are rational functions in the matrix elements of $B_{1},B_{1}',i,j'$.
Similarly the matrix elements of $\tB^{(2,1)}$ are rational functions in the ones of $B_{1},B_{1}',i',j$ defined over the locus where $B_{1},B_{1}'$ do not have a common eigenvalue.

\begin{prop}\label{prop: tensor ADHM}
$($\cite[Theorem 2.15]{Choy2}$)$
The rational map $\sT$ induces the tensor product map $\cM^{\SU(N)}_{n}\times\cM^{\SU(N')}_{n'}\to \cM^{\SU(N+N')}_{nN'+n'N}$.
\end{prop}

We extend $\sT$ to the locus where $i\otimes j'=j\otimes i'=0$ by assigning $\tB^{(1,2)}=\tB^{(2,1)}=0$.
In particular if one of $V,V',W,W'$ is zero, $\sT(x,x')$ is defined.
Note that the extended map is not a morphism, but only a set-theoretic map.

We define the dual of $x$ by
	$$
	x\dual:= (B_{1}\dual,B_{2}\dual,-j\dual,i\dual)\in \bM_{V\dual,W\dual}.
	$$
Let $F$ be the associated monad to $x$ (= a complex of vector bundles on $\pp^2$ via the monad construction \cite[Chap.\ 2]{Lecture}).
Let $\cD^{b}(\pp^{2})$ be the bounded derived category of coherent sheaves on $\pp^{2}$.
It is well-known that the derived dual $F\dual$ in $\cD^b(\pp^2)$ coincides with the monad associated to $x\dual$.

The following two lemmas come from direct calculation of the duals and tensor products:
\begin{lem}\label{lem: stable costable}
$x$ is stable (resp.\ costable) if and only if $x\dual$ is costable (resp.\ stable).
\end{lem}

\begin{lem}\label{lem: tensor product of dual}
If $B_{1},B_{1}'$ have no common eigenvalue,
	$$
	\sT(x\dual,x'\dual)=\sT(x,x')\dual.
	$$
\end{lem}

\begin{cor}\label{cor: ADHM of dual}
For framed vector bundles $F,F'$, the ADHM datum of $(F\otimes F')\dual$ is $\sT(x\dual,x'\dual)$.
 \end{cor}

\proof
The ADHM datum of $(F\otimes F')\dual$ is $\sT(x,x')\dual=\sT(x\dual,x'\dual)$ by Lemma  \ref{lem: tensor product of dual}.
\qed\vskip.3cm

We end up with a lemma, which will be not used in the sequel.

\begin{lem}\label{lem: costable}
Suppose $B_1,B_1'$ have no common eigenvalue.
If $x,x'$ are stable (resp.\ costable), so is $\sT(x,x')$.
\end{lem}

\proof
It is enough to prove only the statement for costability thanks to the above two lemmas.
We use the notation $\sT(x,x')=(\tB_1,\tB_2,\ti,\tj)$.
Let $K$ be a $\tB_1,\tB_2$-invariant subspace in $\Ker \tj$.
By $\tB_1$-invariance we have decomposition
    $$
    K=(K\cap (V_1\otimes W_2))\oplus (K\cap (W_1\otimes V_2)).
    $$
For, each generalized eigenspace of $K$ with respect to $\tB_1$ is a subspace of either $V_1\otimes W_2$ or $W_1\otimes V_2$ due to the assumption on the eigenvalues.

Now $K\cap (V_1\otimes W_2)$ is a $B_1\otimes \Id_{W_2},B_2\otimes\Id_{W_2}$-invariant subspace of $\Ker(j\otimes \Id_{W_2})$.
We claim that $K\cap (V_1\otimes W_2)=0$.
Otherwise there is a nonzero $v\in K\cap (V_1\otimes W_2)$.
Since $\Ker(j\otimes \Id_{W_2})=\Ker(j)\otimes W_2$, we can write $v=\sum_l v_1^{(l)}\otimes w_2^{(l)}$ where $v_1^{(l)}\in \Ker(j)\setminus0$ for each $l$ and $w_2^{(l)}\in W_2$ are linearly independent.
We fix any  $l_0$ among the indices $l$.
Since $x$ is costable there exists a $2$-variable polynomial $f$ such that $f(B_1,B_2)v_1^{(l_0)}$ is not contained in $\Ker(j)$.
Thus $f(\tB_1,\tB_2)v$ has a summand $(f(B_1,B_2)v_1^{(l_0)})\otimes w_2^{(l_0)}$ which does not lie in $\Ker(j\otimes \Id_{W_2})$.
This is contradiction which proves the claim.

Similarly we have $K\cap (W_1\otimes V_2)=0$.
Therefore $K=0$, which shows $\sT(x,x')$ is costable.
\qed\vskip.3cm

 
\subsection{Self-tensor product of ADHM data}\label{subsec: self-tensor}

Now we deduce the self-tensor product of ADHM data.
Thus we need to consider the case when $x=x'$ and the first factor $B_{1}$ has only multiplicity $1$ eigenvalues.
The idea to define $\sT(x,x)\in \bM_{\tV,\tW}$ associated to the self-tensor product framed vector bundle $F^{\otimes2}$ is as follows:
First we consider the ADHM datum $x_{t}:=(B_{1}+t,B_{2},i,j)$ where the first factor is translated by $t\in\cc$.
It is clear that $x_{t}$ is regular and $\sT(x_{t},x)$ is defined for $t\neq0$, but \textit{not} defined for $t=0$.
Let $F_{t}$ denote the associated framed vector bundle to $x_{t}$ so that $\sT(x_{t},x)$ corresponds to $F_{t}\otimes F$ ($t\neq0$) by Proposition \ref{prop: tensor ADHM}.
Next we find a family $\Phi(t)$ in $\GL(\tV)$ parametrized by $\cc^{*}$ such that there exists the limit regular ADHM datum $\lim_{t\to 0}\Phi(t).\sT(x_{t},x)\in \bM_{\tV,\tW}$.
Our main result Proposition \ref{prop: limit} asserts that this limit corresponds to $F^{\otimes 2}$ as framed vector bundles.

We assume $V=\cc$ first for simplicity.
Thus we have $\tV=W^{\oplus2}$ and $\tW=W^{\otimes2}$.
Let $b_{1},b_{2}$ denote the scalars $B_{1},B_{2}$ respectively.
Let $f:=j(1)\in W$ and $e\in W$ be a nonzero vector with $i(e)=1$.
Note that $i$ is nonzero by stability and that $e,f$ are linearly independent since $i(f)=ij(1)=0$.
We decompose $W=\cc\langle e,f\rangle\oplus W_{0}$ where $W_{0}\subset \Ker(i)$.
This gives identification 
	$$
	\begin{aligned}
	W^{\otimes2}=
	&
	\,
	\cc\langle e\otimes e,f\otimes f,e\otimes f,f\otimes e\rangle\oplus (e\otimes W_{0})
	\oplus  (W_{0}\otimes e)
	\\
	&
	\oplus (f\otimes W_{0})\oplus (W_{0}\otimes f)\oplus W_{0}^{\otimes 2}.
	\end{aligned}
	$$
We define a $\cc^{*}$-action $\Phi(t)$ by declaring the weight $\pm1,0$ subspaces as
	$$
	\tV_{-1}=\cc (e,-e),\quad
	\tV_{1}=\cc (f,-f),\quad
	\tV_{0}=\cc\langle (e,e),(f,f)\rangle
	\oplus W_{0}^{\oplus2}.
	$$
We write $\sT(x_{t},x)=(\tB_{1}(t),\tB_{2}(t),\ti,\tj)$.
Note that
	$$
	\tB_{1}(t)= \begin{pmatrix} b_{1}+t & 0 \\ 0 & b_{1} \end{pmatrix},\quad
	\tB_{2}(t)= \begin{pmatrix} b_{2} & t\inv i\otimes j \\ -t\inv j\otimes i & b_{2} \end{pmatrix}
	$$
where all the block matrices are endomorphisms in $\End(W)$.
A direct calculation shows
	\begin{equation}\label{eq: Phi action}
	\begin{array}{lll}
	&
	\Phi(t).\tB_{1}-b_{1}\colon
	&
	(e,-e) \mapsto
	\frac t2 (e,-e)
	+ \frac{t^{2}}2 (e,e),\quad
	(f,-f)\mapsto
	\frac t2 (f,-f)+ \frac12 (f,f),
	\\
	& &
	(e,e) \mapsto
	\frac 12(e,-e) 	+ \frac{t}2 (e,e), \quad
	(f,f) \mapsto
	\frac {t^{2}}2 (f,-f) + \frac t2 (f,f),
	\\
	&
	\Phi(t).\tB_{2}-b_{2}\colon
	&
	(e,-e)\mapsto -(f,f), \quad
	(f,-f)\mapsto 0,
	\\
	& &
	(e,e)\mapsto (f,-f), \quad
	(f,f)\mapsto 0,
	\\
	&
	\Phi(t)\ti\colon
	&
	e\otimes e \mapsto (e,e), \quad
	f\otimes f  \mapsto 0,
	\\
	& &
	e\otimes f  \mapsto
	\frac t2 (f,-f)+\frac 12(f,f), \quad
	f\otimes e \mapsto -\frac t2(f,-f)+\frac 12(f,f),
	\\
	&
	\tj\Phi(t)\inv\colon
	&
	(e,-e)\mapsto t(f\otimes e-e\otimes f), \quad
	(f,-f) \mapsto 0,
	\\
	& &
	(e,e)\mapsto f\otimes e+e\otimes f, \quad
	(f,f) \mapsto 2f\otimes f.
	\end{array}
	\end{equation}
The images of the other basis elements do not change under the $\Phi(t)$-action.
In particular $\Phi(t).\tB_{2}|_{(W_{0})^{\oplus2}}=0$.
As a result $\lim_{t\to 0} \Phi(t).\sT(x_{t},x)$ exists.
We denote the limit ADHM datum by $((\tB_{1})_{0},(\tB_{2})_{0},\ti_{0},\tj_{0})$.

In the general case $k=\dim V\ge1$, we identify $V=\cc^{k}$ using the eigenspace decomposition, hence $\tV=(W^{\oplus2})^{\oplus k}$.
Now we define a $\cc^{*}$-action $\Phi(t)$ on $\tV$ as the diagonal $\cc^{*}$-action on each summand $W^{\oplus2}$ of $\tV$ as before.
More precisely, let $p_{l}\colon V\to \cc,\ q_{l}\colon \cc\to V$ be the $i^{\mathrm{th}}$ projection and the inclusion as the $i^{\mathrm{th}}$ summand respectively.
By the factorization property, $x$ is the $\sigma$-image of the $k$-tuple $(B_{1}^{l},B_{2}^{l},i_{l},j_{l})_{l=1,2,...,k}$ where $(B_{1}^{l},B_{2}^{l},i_{l},j_{l}):=(p_{l}B_{1}q_{l},p_{l}B_{2}q_{l},p_{l}i,jq_{l})$ an ADHM datum in $\bM_{\cc,W}$.
Instead of $e,f$, we use $e_{1},e_{2},...,e_{k}$ and $f_{1},f_{2},...,f_{k}$ in $W$ using $(B_{1}^{l},B_{2}^{l},i_{l},j_{l})$.
Instead of $W_{0}$, we use $W_{0}^{l}\subset \Ker(i_{l})\subset W$ complementary to $\cc\langle e_{l},f_{l}\rangle$ for each $l$.
Then $\Phi(t)$ is set to be the $\cc^{*}$-action on $\tV$ with the weight spaces 
	$$
	\begin{aligned}
	&
	\tV_{-1}=\cc\langle (e_{1},-e_{1}),(e_{2},-e_{2}),...,(e_{k}, -e_{k})\rangle,
	\\
	&
	\tV_{1}=\cc\langle (f_{1}, -f_{1}),(f_{2},-f_{2}),...,
	(f_{k},-f_{k})\rangle,
	\\
	&
	\tV_{0}=
	\bigoplus_{l=1}^{k}\cc\langle (e_{l},e_{l}),(f_{l},f_{l})\rangle\oplus (W_{0}^{l})^{\oplus2}
	\end{aligned}
	$$
Due to the above calculation \eqref{eq: Phi action} applied to each summand $W^{\oplus2}$ of $\tV$, we also obtain the well-defined $t\to0$ limit ADHM datum $\sT(x,x):=((\tB_{1})_{0},(\tB_{2})_{0},\ti_{0},\tj_{0})$.

\begin{prop}\label{prop: limit}
The ADHM datum of the framed vector bundle $F^{\otimes 2}$ is the above $((\tB_{1})_{0},(\tB_{2})_{0},\ti_{0},\tj_{0})$.
\end{prop}

\proof
By \eqref{eq: Phi action} there is no nonzero $(\tB_{1})_{0},(\tB_{2})_{0}$-invariant subspace in
	$$
	\Ker(\tj_{0})
	=\bigoplus_{l=1}^{k}
	\cc\langle(e_{l}, -e_{l}),(f_{l},-f_{l})\rangle.
	$$
Thus $((\tB_{1})_{0},(\tB_{2})_{0},\ti_{0},\tj_{0})$ is costable.

By  \eqref{eq: Phi action} any $(\tB_{1})_{0},(\tB_{2})_{0}$-invariant subspace containing
	$$
	\Image(\ti_{0})
	=\bigoplus_{l=1}^{k}\cc\langle(e_{l},e_{l}),
 	(f_{l},f_{l})\rangle\oplus (W_{0}^{l})^{\oplus2}
	$$
should be $\tV=(W^{\oplus2})^{\oplus k}$.
Thus $((\tB_{1})_{0},(\tB_{2})_{0},\ti_{0},\tj_{0})$ is stable.
Hence it is regular and corresponds to a framed vector bundle.
Since this framed bundle is isomorphic to a limit of $F_{t}\otimes F$, it is isomorphic to $F\otimes F$ as framed bundles.
\qed\vskip.3cm

In general an ADHM datum in $\bM_{V,W}$ of a given framed vector bundle is unique up to $\GL(V)$-action, and the $\cc^{*}$-action $\Phi(t)$ is absorbed in $\GL(\tV)$-action.

We recall that if $F$ has a symplectic structure,  $V,W$ are orthogonal and symplectic vector spaces respectively and $B_{1}=B_{1}^{*},B_{2}=B_{2}^{*},j=i^{*}$.
Thus $\tV,\tW$ are naturally symplectic and orthogonal vector spaces.
We identify $V=\cc^{k}$ the standard orthogonal vector space using the eigenspace decomposition.
Note that $(e_{l},f_{l})_{W}=1$ since $j=i^{*}$.
We set $W_{0}^{l}:=\cc\langle e_{l},f_{l}\rangle^{\perp}$ the orthogonal complement in $W$.
Then the above $\cc^{*}$-action $\Phi(t)$ preserves the symplectic form of $\tV$.
Thus the limit ADHM datum $((\tB_{1})_{0},(\tB_{2})_{0},\ti_{0},\tj_{0})$ is contained in $\bN_{\tV,\tW}$.
Hence Proposition \ref{prop: limit} can be adapted to the symplectic version:

\begin{prop}\label{prop: limit Sp}
Suppose further $F$ is a symplectic bundle.
The ADHM datum of the framed orthogonal bundle $F^{\otimes 2}$ is the above $((\tB_{1})_{0},(\tB_{2})_{0},\ti_{0},\tj_{0})$.
\end{prop}


\subsection{Symmetric product and exterior product of ADHM data}\label{subsec: symm ext ADHM}

In this subsection
we shall find the ADHM datum of the second symmetric product $S^{2}F$ and the second exterior product $\Lambda^{2}F$ of a framed vector bundle $F$.
We also use the assumption that the first factor $B_{1}$ of the ADHM datum $x$ corresponding to $F$ has only multiplicity $1$ eigenvalues.

We consider mutually complementary subspaces $S^{2}W$, $\Lambda^{2}W$ in $W^{\otimes2}$ spanned by the vectors $w\otimes w'+w'\otimes w$, $w\otimes w'-w'\otimes w$ respectively.
Let
	$$
	V_{S}:=\sum_{h}h((\tB_{1})_{0},(\tB_{2})_{0})\ti_{0}(S^{2}W),\quad
	V_{E}:=\sum_{h}h((\tB_{1})_{0},(\tB_{2})_{0})\ti_{0}(\Lambda^{2}W)
	$$
where $h$ runs over the $2$-variable polynomials.
We denote the involutive and anti-involutive subspaces in $W^{\oplus2}$ by
	$$
	\Delta^{\pm} W:=\left\{(w,\pm w)\in W^{\oplus2}\right\}.
	$$
We define similar subspaces in $W_{0}^{\oplus 2}$ using the notation $\Delta^{\pm}$.
\begin{lem}\label{lem: VS VE}
$V_{S},V_{E}$ are identified as
	$$
	V_{S}=\bigoplus_{l=1}^{k}\cc\langle e_{l},f_{l}\rangle^{\oplus2} \oplus \Delta^{+} W_{0}^{l},\quad
	V_{E}=\bigoplus_{l=1}^{k} \Delta^{-} W_{0}^{l}.
	$$
\end{lem}
\proof
We check the first identification.
For each $l$ there is decomposition
	$$
	\begin{aligned}
	S^{2}W=
	&\,
	\cc\langle e_{l}\otimes e_{l},f_{l}\otimes f_{l},e_{l}\otimes f_{l}+f_{l}\otimes e_{l}\rangle
	\oplus \left\{e_{l}\otimes w+ w\otimes e_{l}\middle| w\in W_{0}^{l}\right\}
	\\
	&
	\oplus \left\{f_{l}\otimes w+ w\otimes f_{l}\middle| w\in W_{0}^{l}\right\}
	\oplus S^{2}W_{0}^{l}
	\end{aligned}
	$$
We denote by $(S^{2}W)_{m}$, $m=1,2,3,4$ the above summands in order.
Let $p_{l}\colon \tV=(W^{\oplus2})^{\oplus k}\to W^{\oplus2}$ be the $l^{\mathrm{th}}$ projection.

Via $p_{l}$, the sum $\sum_{h}h((\tB_{1})_{0},(\tB_{2})_{0})\ti_{0}((S^{2}W)_{1})$ projects onto $\cc\langle e_{l}\otimes e_{l},f_{l}\otimes f_{l},e_{l}\otimes f_{l},f_{l}\otimes e_{l}\rangle$.
The sum for $(S^{2}W)_{2}$ projects onto $\Delta^{+}W_{0}^{l}$ since it coincides with the set of $p_{l}\ti(e_{l}\otimes w+w\otimes e_{l})=(w,w),\ w\in W_{0}$.
Both sums for $(S^{2}W)_{3},(S^{2}W)_{4}$ project to zero.
Therefore we get 
	$$
	p_{l}(V_{S})=\cc\langle e_{l},f_{l}\rangle\oplus \Delta^{+}W_{0}^{l}.
	$$ 
We use $V_{S}=\bigoplus_{l=1}^{k}p_{l}(V_{S})$, which comes from that the $l^{\mathrm{th}}$ summands $W^{\oplus2}$ of $\tV$ form the eigenspaces of $\tB_{1}$ with distinct eigenvalues respectively.
This identifies $V_{S}$ as in the statement. 

The identification of $V_{E}$ is similarly checked using decomposition
	$$
	\begin{aligned}
	\Lambda^{2}W
	 =
	 & \, \cc\langle e_{l}\otimes f_{l}-f_{l}\otimes e_{l}\rangle\oplus
	\left\{e_{l}\otimes w- w\otimes e_{l}\middle| w\in W_{0}^{l}\right\}
	\\
	&
	\oplus \left\{f_{l}\otimes w- w\otimes f_{l}\middle| w\in W_{0}^{l}\right\}
	\oplus \Lambda^{2}W_{0}^{l}.
	\end{aligned}
	$$
Via $p_{l}$, the sum $\sum_{h}$ as above for the first summand projects to zero by \eqref{eq: Phi action}.
The sum for the second summand projects onto $\Delta^{-}W_{0}$. 
Both sums for the other two summands project to zero.
This identifies $V_{E}$ as in the statement.
\qed\vskip.3cm

 \begin{rk}\label{rk: second Chern}
$\dim V_{S}=c_{2}(S^{2}F)=(\rank F+2)c_{2}(F)$ and $\dim V_{E}=c_{2}(\Lambda^{2}F)=(\rank F-2)c_{2}(F)$.
\end{rk}
\begin{thm}\label{th: symm ext}
The restrictions of $((\tB_{1})_{0},(\tB_{2})_{0},\ti_{0},\tj_{0})$ to $(V_{S},S^{2}W)$, $(V_{E},\Lambda^{2}W)$
are the ADHM data of $S^{2}F,\Lambda^{2}F$ respectively.
\end{thm}
\proof
We first need to check the restrictions are well-defined, i.e., $\tj_{0}(V_{S})\subset S^{2}W,\ \tj_{0}(V_{E})\subset \Lambda^{2}W$.
By \eqref{eq: Phi action}, for any $l$ we have $\tj_{0}(\cc\langle e_{l},f_{l}\rangle^{\oplus2})=\cc\langle f_{l}\otimes e_{l}+e_{l}\otimes f_{l},f_{l}\otimes f_{l}\rangle$.
We also have
	$$
	\tj_{0}(\Delta^{+}W_{0}^{l})=
	\tj(\Delta^{+}W_{0}^{l})=\left\{f_{l}\otimes w+w\otimes f_{l}\middle| w\in W_{0}\right\}.
	$$
Thus $\tj_{0}(V_{S})\subset S^{2}W$ by Lemma \ref{lem: VS VE}.
Similarly
	\begin{equation}\label{eq: j0}
	\tj_{0}(\Delta^{-}W_{0}^{l})
	=\tj(\Delta^{-}W_{0}^{l})
	=\left\{f_{l}\otimes w-w\otimes f_{l}\middle| w\in W_{0}\right\}\subset \Lambda^{2}W.
	\end{equation}
Thus $\tj_{0}(V_{E})\subset \Lambda^{2}W$.

Now $((\tB_{1})_{0},(\tB_{2})_{0},\ti_{0},\tj_{0})$ is the direct sum of the two restrictions  to $(V_{S},S^{2}W)$, $(V_{E},\Lambda^{2}W)$ as quiver representations.
Since $((\tB_{1})_{0},(\tB_{2})_{0},\ti_{0},\tj_{0})$ is a regular element (Proposition \ref{prop: limit}), these restrictions are also regular.
Let $F_{S},F_{E}$ be the corresponding framed vector bundles.
We have $F^{\otimes2}=F_{S}\oplus F_{E}$.
On the other hand there is an inclusion $F_{S}\subset S^{2}F$ as framed sheaves, because the frames of $F_{S},S^{2}F$ coincide and $F_{S}$ is the minimal framed subsheaf in $F^{\otimes2}$ with such a frame due to Lemma  \ref{lem: VS VE}.
Similarly $F_{E}\subset \Lambda^{2}F$.
By the decomposition $F^{\otimes2}=S^{2}F\oplus\Lambda^{2}F$ these inclusions are nothing but $F_{S}=S^{2}F,\ F_{E}=\Lambda^{2}F$.
\qed\vskip.3cm


\subsection{Construction of the tensor, symmetric and exterior product morphisms in terms of ADHM data}\label{subsec: construction of isom}

We construct the morphisms $F\mapsto \ad F,(\Lambda^{2}F)_{0},\Lambda^{2}F$ in terms of ADHM data $x=(B_{1},B_{2},i,j)$ of $F$ introduced in the beginning of this section.
To be precise we will give the explicit forms of the morphisms only over the locus 
	$$
	\mu\inv(0)^{\reg}_{0}:=\left\{x\in \mu\inv(0)^{\reg}\middle|
	\mbox{$B_{1}$ has only the multiplicity $1$ eigenvalues}\right\}.
	$$
These morphisms extend to the morphisms between the framed vector bundles due to the factorization property since the quasi-affine GIT quotient of $\mu\inv(0)^{\reg}_{0}$  is Zariski dense open in $\cM^{K}_{n}$.

We denote the moment map on $\bM_{\tV,\tW}$ or $\bN_{\tV,\tW}$ by $\tmu$.

\begin{prop}\label{prop: self tensor morphism}
$(1)$
The morphism $\cM^{\SU(N)}_{n}\to \cM^{\SU(N^{2})}_{2nN},\ F\mapsto F^{\otimes2}$ is induced by the morphism $\mu\inv(0)^{\reg}_{0}\to \tmu\inv(0)^{\reg}/\GL(\tV),\ x\mapsto [\sT(x,x)]$.

$(2)$
The morphism $\cM^{\USp(N/2)}_{n}\to \cM^{\SO(N^{2},\rr)}_{2nN},\ F\mapsto F^{\otimes2}$ is induced by the morphism $\mu\inv(0)^{\reg}_{0}\to \tmu\inv(0)^{\reg}/\Sp(\tV),\ x\mapsto [\sT(x,x)]$.
\end{prop}

\proof
It suffices to observe that $\Phi(t)$ is absorbed in $\GL(\tV)$ or $\Sp(\tV)$ by construction (Proposition \ref{prop: limit Sp}).
\qed\vskip.3cm

We construct first the morphisms $F\mapsto S^{2}F,\Lambda^{2}F$ in terms of ADHM data.
We also use the same notation $\tmu$ for the moment map on $\bM_{V_{S},S^{2}W}$ or $\bM_{V_{E},\Lambda^{2}W}$.
We fix a triple $(x_{\flat},V_{S}^{\flat},V_{E}^{\flat})$ where $x_{\flat}=(B_{1}^{\flat},B_{2}^{\flat},i_{\flat},j_{\flat})\in \mu\inv(0)^{\reg}\subset \bM_{V,W}$ and  $V_{S}^{\flat},V_{E}^{\flat}$ are subspaces of $\tV$ constructed from $x_{\flat}$ and choices of $e_{l},f_{l}$ as before.
For each triple $(x,V_{S},V_{E})$ we choose any $\phi\in\GL(\tV)$ with $V_{S}^{\flat}= \phi(V_{S}),\ V_{E}^{\flat}=\phi(V_{E})$.
The restrictions of $\phi.\sT(x,x)$ to $(V_{S}^{\flat},S^{2}W)$ and $(V_{E}^{\flat},\Lambda^{2}W)$ are elements of $\bM_{V_{S}^{\flat},S^{2}W}$ and $\bM_{V_{E}^{\flat},\Lambda^{2}W}$ respectively.
Their classes in the $\GL(V_{S}^{\flat})$- and $\GL(V_{E}^{\flat})$-quotients $\bM_{V_{S}^{\flat},S^{2}W}\git \GL(V_{S}^{\flat})$ and $\bM_{V_{E}^{\flat},\Lambda^{2}W}\git \GL(V_{E}^{\flat})$ do not depend on the choice of $\phi$.
This induces morphisms between instanton spaces:

\begin{prop}\label{prop: symm ext morphism}
$(1)$
The symmetric product morphism $\cM^{\SU(N)}_{n}\to \cM^{\SU(N')}_{n'}$, $F\mapsto S^{2}F$ is induced by the morphism $\mu\inv(0)^{\reg}_{0}\to \tmu\inv(0)^{\reg}/\GL(V_{S}^{\flat})$, $x\mapsto [\phi.\sT(x,x)|_{(V_{S}^{\flat},S^{2}W)}]$ where $N'=\frac{N(N+1)}2$ and $n'=n(N+2)$.

$(2)$
The exterior product morphism $\cM^{\SU(N)}_{n}\to \cM^{\SU(N')}_{n'}$, $F\mapsto \Lambda^{2}F$ is induced by the morphism $\mu\inv(0)^{\reg}_{0}\to\tmu\inv(0)^{\reg}/\GL(V_{E}^{\flat})$, $x\mapsto[\phi.\sT(x,x)|_{(V_{E}^{\flat},\Lambda^{2}W)}]$ where $N'=\frac{N(N-1)}2$ and $n'=n(N-2)$.
\end{prop}

We get back to the construction of morphisms $F\mapsto \ad F,(\Lambda^{2}F)_{0},\Lambda^{2}F$ in terms of ADHM data.

\begin{thm}\label{th: isom in ADHM}
$(1)$
The isomorphism $\cM^{\SU(2)}_{n}\cong \cM^{\SO(3,\rr)}_{n}$, $F\mapsto \ad F$ is given by the symmetric product morphism in Proposition \ref{prop: symm ext morphism} (1).

$(2)$
The isomorphism $\cM^{\USp(2)}_{n}\cong \cM^{\SO(5,\rr)}_{n}$ is given by $x\mapsto \phi.\sT(x,x)|_{(V_{E}^{\flat},\Ker(\omega))}$ where $\omega\colon \Lambda^{2}W\to\cc$ is the symplectic form on $W$.

$(3)$
The isomorphism $\cM^{\SU(4)}_{n}\cong \cM^{\SO(6,\rr)}_{n}$ is given by $x\mapsto \phi.\sT(x,x)|_{(V_{E}^{\flat},\Lambda^{2}W))}$.
\end{thm}

In the statements (2) and (3) we need to explain additional structures on $V_{E}^{\flat},\Lambda^{2}W,\phi$, etc.
In the item (2), $V,W$ are orthogonal and symplectic vector spaces respectively.
Recall that $\tV,\tW$ have the induced symplectic and orthogonal structures respectively.
Recall also that for a given $x\in \mu\inv(0)^{\reg}_{0}$, further choice of $e_{l}$ defines $V_{S},V_{E}$.

\begin{lem}\label{lem: symplectic str of tV}
$V_{S},V_{E}$ are symplectic subspaces of $\tV$.
\end{lem}

\proof
We check $V_{E}$ is a symplectic subspace of $\tV$.
We identify $V=\cc^{k}$ the standard orthogonal vector space and $\tV=(W\otimes \cc^{2})^{\oplus k}$.
For each $l$, $W_{0}^{l}$ is a symplectic subspace in $W\otimes \cc^{2}$ by the definition.
Thus so is $V_{E}=\bigoplus_{l}\Delta^{-}W_{0}^{l}$ in $\tV$.

One can check similarly that $V_{S}$ is a symplectic subspace.
\qed\vskip.3cm

The orthogonal structure on $\tW=W^{\otimes 2}$ is given by the induced isomorphism $W^{\otimes 2}\cong (W\dual)^{\otimes 2}=(W^{\otimes 2})\dual$.
The subspace $\Lambda^{2}W$ is an orthogonal subspace in $\tW$.
This orthogonal structure of $\Lambda^{2}W$ is given by the $4$-form $\wedge^{2}\omega$.
$\Ker(\omega)$ is an orthogonal subspace of $\Lambda^{2}W$.
In (2) we need the restriction $\sT(x,x)$ to $(V_{E},\Ker(\omega))$ is well-defined.
This amounts to $\tj_{0}(V_{E})\subset \Ker(\omega)$, which follows immediately from \eqref{eq: j0}.
We set $\phi$ to be any symplectic isomorphism of $\tV$ satisfying $V_{S}^{\flat}=\phi(V_{S}),V_{E}^{\flat}=\phi(V_{E})$.

In the case (3) we give any orthogonal structure on $V$ and then identify $V=\cc^{k}$ the standard orthogonal vector space.
For each ADHM datum $x\in \mu\inv(0)_{0}$ we use the notation $f_{1}^{x},f_{2}^{x},...,f_{k}^{x}\in W$ for the $j$-images of the standard basis elements of $V$ to emphasize $x$.
We choose a generic symplectic structure $\omega^{x}\colon \Lambda^{2}W\to\cc$ such that there exist $e_{1}^{x},e_{2}^{x},...,e_{k}^{x}\in W$ satisfying  $\omega^{x}(e_{l}^{x},f_{l}^{x})=1$.
The existence of such $\omega^{x}$ is left as an exercise.
Note that $\omega^{x}(e_{l}^{x},f_{l}^{x})=1$ is no more automatic because $j$ is not necessarily $i^{*}$ with respect to $\omega^{x}$.
We further impose a condition on $\omega^{x}$: the $4$-forms $\wedge^{2}\omega^{x}$ do not depend on $x$.
This condition is satisfied simply by scalar multiplications of $\omega^{x}$ and $e_{l}^{x}$ since $\Lambda^{4}W$ is $1$-dimensional.
Thus $\Lambda^{2}W$ attains the orthogonal structure independent of $x$.
It is now easy to check that $\sT(x,x)$ restricts to an element of $\bN_{V_{E},\Lambda^{2}W}$.
We choose $\phi\in \Sp(\tV)$ with $V_{S}^{\flat}=\phi(V_{S}),V_{E}^{\flat}=\phi(V_{E})$ as before.
Hence we have the morphism in the statement $\mu\inv(0)_{0}\to \bN_{V_{E}^{\flat},\Lambda^{2}W}\git \Sp(V_{E}^{\flat})$.

The morphisms in Theorem \ref{th: isom in ADHM} are now well-defined.
We are ready to prove the theorem. 
\vskip.3cm

\textit{Proof of Theorem \ref{th: isom in ADHM}.}
(1) 
Since $\SU(2)=\USp(1)$, there is an isomorphism $F\cong F\dual$ corresponding to the symplectic structure of $F$.
Thus we have the isomorphism $\End(F)\cong F^{\otimes2}$ and this restricts to $\ad F\cong S^{2}F$.
Hence the map $F\mapsto \ad F$ is given by Proposition \ref{prop: symm ext morphism} (1).

(2), (3)
The proof of the item (3) is a first half of that of (2), so we prove only (2).
By similar arguments in the proofs of Proposition \ref{prop: self tensor morphism} (2) and Proposition \ref{prop: symm ext morphism} (2), we obtain a morphism $\mu\inv(0)^{\reg}_{0}\to \tmu\inv(0)^{\reg}/\Sp(V_{E})$ which induces $\cM^{\USp(2)}_{n}\to \cM^{\SO(6,\rr)}_{2n}$, $F\mapsto \Lambda^{2}F$.
The rest follows from the fact that $(\Lambda^{2}F)_{0}=\Ker(\Lambda^{2}F\to \cO)$ has the ADHM datum $\sT(x,x)|_{(V_{E},\Ker(\omega))}$.
To check this fact, we notice that the morphism between the ADHM data induced by $V_{E}\oplus\Lambda^{2}W \stackrel{0\oplus \omega}\to 0\oplus\cc$ becomes the morphism  $\Lambda^{2}F\to \cO$.
\qed\vskip.3cm


\appendix


\section{Character of $\cc[\rho\inv(0)]^{\Sp(1)}$}\label{subsec: character}

In this section we compute the Hilbert series of $\rho\inv(0)\git \Sp(1)$ where $\rho\colon\Hom(\cc^3,\cc^2)\to \fsp(1),\ i\mapsto ii^{*}$ and $\cc^3,\cc^2$ are the standard orthogonal and symplectic vector spaces respectively.

Since $\fo(3)$ is spanned by the weight $-1,0,1$ vectors with respect to a maximal torus $\SO(3)$, we denote by $x,y,z$ the dual basis elements.
Thus we identify $\cc[\fo(3)]=\cc[x,y,z]$.
By the first fundamental theorem of invariant theory $\Hom(\cc^{3},\cc^{2})\git \Sp(1)$ is $\SO(3)$-equivariantly isomorphic to $\fo(3)$.
Hence $\cc[\rho\inv(0)\git \Sp(1)]$ is $\SO(3)$-equivariantly isomorphic to a quotient algebra of $\cc[\fo(3)]$.
By the above computation of Hilbert series we will deduce the following:

\begin{cor}\label{cor: Hilbert series}
There is an $\SO(3)$-equivariant isomorphism 
	$$
	\cc[\rho\inv(0)]^{\Sp(1)}\cong\cc[x,y,z]/(x,y,z)^{2}.
	$$ 
\end{cor}

The proof will appear in \S\ref{subsec: proof of corollary} after identifying the $\SO(3)$-representation $\cc[\rho\inv(0)]^{\Sp(1)}$ as a torus character.


\subsection{Torus character of $\cc[\rho\inv(0)]^{\Sp(1)}$}

For an algebraic group $G$ we denote by $R(G)$ the ring of isomorphism classes of finite dimensional $G$-representations.
For a given $G$-representation $V$, its $G$-character is denoted by $\cX_{V}\colon G\to \cc,\ g\mapsto \tr(g\colon V\to V)$.
Thus we have character map $\cX_{\bullet}\colon R(G)\to \cc[G]$.
If $G$ is a reductive group, let $T$ be a maximal torus.
The composite of the character map $\cX_{\bullet}\colon R(G)\to \cc[G]$ and the restriction $\cc[G]\to \cc[T]$ is known to be injective.
And its image is the Weyl group-invariant ring $\cc[T]^{W}$. 
We denote this composite by the same symbol $\cX$ and call torus character.

We want to add some infinitely dimensional representations to $R(G)$, e.g.\  coordinate ring of a scheme with nonzero dimension.
We narrowly focus on the example $\cc[\rho\inv(0)]$.
We set $G=\Sp(1)\times\SO(3)$ from now on.
Let $T=T_{\Sp(1)}\times T_{\SO(3)}$ where $T_{\Sp(1)}$ and $T_{\SO(3)}$ are maximal tori of $\Sp(1)$ and $\SO(3)$ respectively.
Regarding $\cc^{3}$ and $\cc^{2}$ as the vector representations of $\SO(3)$ and $\Sp(1)$ respectively and then defining a $\cc^{*}$-action with only weight $1$, $\Hom(\cc^{3},\cc^{2})$ is a $G\times \cc^{*}$-representation.  
So is its coordinate ring $\cc[\Hom(\cc^{3},\cc^{2})]=\Sym(\Hom(\cc^{3},\cc^{2})\dual)$ the total symmetric product.
The torus character $\cX_{\Hom(\cc^{3},\cc^{2})}=(1+t+t\inv)(z+z\inv)q$.
Here $t,z,q$ are torus characters corresponding to the $1$-dimensional representations of $T_{\SO(3)},T_{\Sp(1)},\cc^{*}$ with weight $1$ after identifying $T_{\SO(3)}\cong \cc^{*},\ T_{\Sp(1)}\cong \cc^{*}$.
Note that $\cc[T]=\zz[t^{\pm1},z^{\pm1},q^{\pm1}]$.

Each $T$-weight space of $\Sym(\Hom(\cc^{3},\cc^{2})\dual)$ is finite dimensional.
Hence $\cX_{\cc[\Hom(\cc^{3},\cc^{2})]}$ is an element in both completed rings
	$$
	\hR(T\times T'):=R(T\times T')[[q\inv]],\quad \hR(G\times T'):=R(G\times T')[[q\inv]].
	$$
Here and hereafter we identified $z,t,q$ with their corresponding representations.
In $\hR(T\times T')$ we have
	$$
	\cX_{\cc[\Hom(\cc^{3},\cc^{2})]}=\frac{1}{P(zq\inv)P(z\inv q\inv)}.
	$$
where $P(x):=(1-x)(1-tx)(1-t\inv x)$ and the RHS is understood as the formal series expansion in $q\inv$.

We regard the adjoint $\Sp(1)$-representation $\fsp(1)$ as a $G\times T'$-representation with the trivial $\SO(3)$-action and the weight $2$ $T'$-action.
Thus $\rho$ becomes $G\times T'$-equivariant.
By finite dimensionality of $T'$-weight spaces, $\cX_{\cc[\rho\inv(0)]}$ is an element in both $\hR(T\times T'),\hR(G\times T')$.
The pull-back of $\fsp(1)\dual$ via $\rho$ generates the defining ideal of $\rho\inv(0)$.
These can be seen as sections of the trivial vector bundle $\cV:=\fsp(1)\dual\times \Hom(\cc^{3},\cc^{2})$ over $\Hom(\cc^{3},\cc^{2})$.
Since $\rho\inv(0)$ is a complete intersection (\cite[Theorem 4.1 (2)]{Choy}), the Koszul complex $\Lambda^{\bullet}\cV$ of $\cO_{\Hom(\cc^{3},\cc^{2})}$-modules is equal to $\cO_{\rho\inv(0)}$ as classes of  the $G\times T'$-equivariant Grothendieck group $K^{G\times T'}(\Hom(\cc^{3},\cc^{2}))$.
Hence  we have 
	$$
	\begin{aligned}
	&
	\cX_{\cc[\rho\inv(0)]}=\sum_{l}(-1)^{l}\cX_{H^{0}(\Lambda^{l}\cV)} 
	=\sum_{l}(-1)^{l}\cX_{\Lambda^{l}\fsp(1)\dual}\cdot\cX_{\cc[\Hom(\cc^{3},\cc^{2})]}
	\\
	& 
	=\frac{(1-q^{-2})(1-z^{2}q^{-2})(1-z^{-2}q^{-2})}{P(zq\inv)P(z\inv q\inv)} 
	\end{aligned}
	$$
in $\hR(T\times T')$ where $H^{0}(\Lambda^{l}\cV)$ denotes the space of sections of $\Lambda^{l}\cV$ over $\Hom(\cc^{3},\cc^{2})$.

The character of the invariant subspace $\cc[\rho\inv(0)]^{\Sp(1)}$ becomes an element of a completed ring $\hR(T_{\SO(3)}\times T')=R(T_{\SO(3)}\times T')[[q\inv]]$.
It is computed by Weyl's integration formula as
	$$
	\cX_{\cc[\rho\inv(0)]^{\Sp(1)}}=\frac12 \ointctrclockwise_{|z|=1} \frac{dz}{2\pi\sqrt{-1} z}\cdot N(1)\cdot \frac{(1-q^{-2})N(q^{-2})}{P(zq\inv)P(z\inv q\inv)}
	$$
(the counter-clockwise integral), where $N(x):=(1-xz^{2})(1-xz^{-2})$ the Jacobian of the finite map $T_{\Sp(1)}\times \Sp(1)/T\to \Sp(1),\ (t,gT)\mapsto g\inv t g$.
We notice that the integrand formal series in $q\inv$ converges in the the range $|q|\gg1$.
Thus it is also regarded as a rational function in $|q|\gg1$.
The denominator of this rational function has the zeros $z=0,\frac1{qt},\frac tq$ inside $|z|=1$. 
By direct residue computation we get
	\begin{equation}\label{eq: torus}
	\cX_{\cc[\rho\inv(0)]^{\Sp(1)}}=1+tq^{-2}+q^{-2}+t\inv q^{-2}
	\end{equation}
in $\hR(T_{\SO(3)}\times T')$.
Since this character is a finite sum, it is also defined in $R(T_{\SO(3)}\times T')$.

\begin{rk}\label{rk: rho 0}
By the second fundamental theorem of invariant theory, the reduced scheme $(\rho\inv(0)\git\Sp(1))_{\mathrm{red}}$ coincides with the zero in $\fo(3)$. 
But in the above we checked $\rho\inv(0)\git\Sp(1)$ is not reduced.
\end{rk}

\subsection{Proof of Corollary \ref{cor: Hilbert series}}\label{subsec: proof of corollary}
We give the weight $2$ $T'$-action on $\fo(3)$.
Then we have $\cX_{\cc[\fo(3)]}=\frac1{(1-tq^{-2})(1-q^{-2})(1-t\inv q^{-2})}$ in $\hR(T_{\SO(3)}\times T')$.
In the statement of Corollary \ref{cor: Hilbert series}, the variables $x,y,z$ can be set to be the bases of the $1$-dimensional subrepresentations $tq^{-2},q^{-2},t\inv q^{-2}$ respectively.

By \eqref{eq: torus}, $\cc[\rho\inv(0)]^{\Sp(1)}$ coincides the subspace of $\cc[\fo(3)]$ spanned by $1,x,y,z$.
On the other hand by the first fundamental theorem of invariant theory, $\Hom(\cc^{2},\cc^{3})\git \Sp(1)$ is an $\SO(3)\times T'$-invariant closed subscheme of $\fo(3)$.
As a quotient ring of $\cc[\fo(3)]$, $\cc[\Hom(\cc^2,\cc^3)]^{\Sp(1)}$ should be $\cc[\fo(3)]/(x,y,z)^{2}$.
This proves the corollary.


\end{document}